\documentclass[leqno,12pt]{amsart}
\usepackage{amsmath,amstext,amssymb,amsopn,amsthm,mathrsfs}
\usepackage[top=4cm, bottom=3cm, left=3cm, right=3cm]{geometry}
\usepackage[shortlabels]{enumitem}
\usepackage{amsfonts}

\usepackage[T1]{fontenc}
\usepackage[makeroom]{cancel}
\usepackage{color}
\usepackage[dvips]{graphicx}
\usepackage{graphicx}
\usepackage{ulem}
\graphicspath{ {./images/} }

\usepackage{psfrag}
\usepackage{picinpar}
\usepackage[hyphens]{url}        %paquete para el hypervinculo con el indice, referencias, etc.
\usepackage[breaklinks,colorlinks=true,linkcolor=blue,citecolor=blue, urlcolor=blue]{hyperref}                                                 %configuracion para el hypervinculo
\usepackage[utf8x]{inputenc}
\usepackage{txfonts}

\usepackage{longtable}%paquete para
\newtheorem{thm}{Theorem}[section]
\newtheorem{df}{Definition}[section]
\newtheorem{cor}{Corollary}[section]
\newtheorem{lem}{Lemma}[section]
\newtheorem{prob}{Proposition}[section]

\numberwithin{equation}{section}

\begin{document}

\title[Fourier Series] {NOTES ON HARMONIC ANALYSIS
PART II: \\
THE FOURIER SERIES
}
\author{Kecheng Zhou}
\address{Department of Mathematics, California State University, Sacramento, CA, 95819, USA.}
\email{[Kecheng Zhou]zhouk@csus.edu}
   \author{M. Vali Siadat}
\address{Department of Mathematics and Statistics, Loyola University Chicago,
Chicago, IL, 60660, USA.}
\email{[M. Vali Siadat]msiadat@luc.edu}


\begin{abstract}
Fourier Series is the second of monographs we present on harmonic analysis. Harmonic analysis is one of the most fascinating areas of research in mathematics. Its centrality in the development of many areas of mathematics such as partial differential equations and integration theory and its many and diverse applications in sciences and engineering fields makes it an attractive field of study and research.

The purpose of these notes is to introduce the basic ideas and theorems of the subject to students of mathematics, physics, or engineering sciences. Our goal is to illustrate the topics with utmost clarity and accuracy, readily understandable by the students or interested readers. Rather than providing just the outlines or sketches of the proofs, we have actually provided the complete proofs of all theorems. This approach will illuminate the necessary steps taken and the machinery used to complete each proof.

The prerequisite for understanding the topics presented is the knowledge of Lebesgue measure and integral. This will provide ample mathematical background for an advanced undergraduate or a graduate student in mathematics.
\end{abstract}

\maketitle

\section{Definitions and important results}

\begin{df}
The set of all complex numbers of modulus 1 is denoted by 
$${\bf T}=\{z=e^{ix}:x\in R\}.$$
\end{df}

${\bf T}$ is a compact abelian group with binary operation: complex
multiplication and topology: open arcs $\{e^{ix} : x\in (a,b)\}.$

Define the periodic function $F(x)$ on $R$ by 
\[
F(x)=f(e^{ix}),\quad x\in R,
\]
where $f$ is a function on ${\bf T}.$ 

Let $\chi $ be the identity function on $T$, i.e.,
\[
\chi (e^{ix})=e^{ix},\quad x\in R. 
\]
Clearly, $F(x)=\cos \,x+i\sin \,x$ satisfies $F(x)=\chi (e^{ix})$ for all $
x\in R.$

\begin{df}
\[
{\bf L^p(T)}=\{f\,\text{defined on} \,T: \int |f|^p d\sigma<\infty\}, 
\]
where 
\[
\int |f|^p d\sigma=\int_{-\pi}^{\pi}|f(e^{ix})|^p \frac{dx}{2\pi} =\frac{1}{
2\pi}\int_{-\pi}^\pi |F(x)|^pdx. 
\]
\end{df}

\begin{thm}
\[
{\bf L^{p}(T)\supset L^{r}(T)},\quad \mbox{if $p<r,$}\,\, that\,\, is \quad \
||f||_{p}\leq ||f||_{r}. 
\]
\end{thm}

{\bf Proof:}\quad Using H\"{o}lder's inequality, we have: ($q=r/p>1$) 
\[
\int |f|^{p}d\sigma =\int |f|^{p}\cdot 1d\sigma \leq (\int
(|f|^{p})^{q}d\sigma )^{1/q}(\int 1^{q^{\prime }}d\sigma )^{1/q^{\prime
}}=(\int |f|^{r}d\sigma )^{p/r}<\infty,
\]
where $\frac{1}{q}+\frac{1}{q'}=1.$

\begin{df}
If $f\in L^{1}(T)$, define Fourier coefficients of $f$ as follows: for $
n=0,\pm 1,\pm 2,\cdots ,$ 
\[
a_{n}(f)=\int f\chi ^{-n}d\sigma =\frac{1}{2\pi }\int_{\pi }^{\pi
}F(x)e^{-inx}dx. 
\]
We now formally introduce the series
\[
f(e^{ix})\sim \sum_{n=-\infty }^{\infty }a_{n}(f)e^{inx}. 
\]
The series is called the Fourier series of $f$. Whenever we speak of convergence or summability of a Fourier series, we are
always concerned with the limit, ordinary or generalized, of the 
symmetric partial sums.
\end{df}

\begin{thm}
Let $f$ be a function on $T$ defined as $f=\sum_{-\infty}^{
\infty}a_n\chi^n$ so that the right-hand side series converges uniformly. Then $
a_n(f)=a_n$ for all $n$. That is, the Fourier series of $f$ is the  series $\sum_{-\infty}^{\infty}a_n\chi^n$.
\end{thm}

{\bf Proof:}\quad To compute $a_n(f)$, we integrate $\int f\chi^{-n}=\int
(\sum_k a_k\chi^k)\chi^{-n}.$ To integrate the latter, we integrate
term-by-term. It is worth noting that if $\sum_{-\infty}^{\infty}a_n\chi^n$
converges uniformly on $T$ for some ordering of the series, then $a_n(f)=a_n$
and the series is the Fourier series of $f$. Assume that the series
converges uniformly to $g(x)$ on $T$ for some other ordering. Then $
a_n(g)=a_n$ so that $a_n(f)=a_n(g)$ and $f=g, a.e.$ Since $f$ and $g$ are
continuous, $f=g$ everywhere on $T.$

\begin{thm}
\begin{enumerate}
\item 
\[
|a_n(f)|\leq ||f||_1,\,\, \forall n 
\]

\item A more precise result: (Bessel's inequality) If $f\in L^{2}(T),$ then 
\[
\sum |a_{n}(f)|^{2}\leq ||f||_{2}^{2}. 
\]
\end{enumerate}
\end{thm}

{\bf Proof:}\quad (1) is trivial. As to (2), let 
\[
f(e^{ix})\sim \sum_{n=-\infty }^{\infty }a_{n}(f)e^{inx} 
\]
be the Fourier series of $f$. Define the symmetric partial sums as
\[
f_{N}=\sum_{n=-N}^{N}a_{n}(f)\chi ^{n}. 
\]
We can obtain (2) directly from the following computation: 
\begin{eqnarray*}
||f-f_{N}||_{2}^{2} &=&||f||_{2}^{2}+||f_{N}||_{2}^{2}-2Re\int f\overline{
f_{N}}d\sigma \\
&=&||f||_{2}^{2}+\sum_{-N}^{N}|a_{n}(f)|^{2}-2\sum_{-N}^{N}Re(\overline{
a_{n}(f)}\int f\overline{\chi }^{n}d\sigma ) \\
&=&||f||_{2}^{2}+\sum_{-N}^{N}|a_{n}(f)|^{2}-2\sum_{-N}^{N}Re(\overline{
a_{n}(f)}\int f{\chi }^{-n}d\sigma ) \\
&=&||f||_{2}^{2}+\sum_{-N}^{N}|a_{n}(f)|^{2}-2\sum_{-N}^{N}Re(\overline{
a_{n}(f)}a_{n}(f)) \\
&=&||f||_{2}^{2}-\sum_{-N}^{N}|a_{n}(f)|^{2}.
\end{eqnarray*}
Thus,
\[
||f-f_{N}||_{2}^{2}=||f||_{2}^{2}-\sum_{-N}^{N}|a_{n}(f)|^{2} 
\]
so that 
\[
||f||_{2}^{2}-\sum_{-N}^{N}|a_{n}(f)|^{2}\geq 0,\quad \forall N. 
\]
We have (2) by taking $N\rightarrow \infty .$

\begin{cor}
If $f\in L^{2}(T),$ then the Fourier partial sums $f_{N}\rightarrow f$ in $
L^{2}(T)$ iff 
\[
||f||_{2}^{2}=\sum |a_{n}(f)|^{2}. 
\]
This equality is called the Parseval relation.
\end{cor}

{\bf Proof:}\quad The corollary follows directly from 
\[
||f-f_N||_2^2= ||f||_2^2-\sum_{-N}^N|a_n(f)|^2. 
\]

\begin{thm}
(Riesz-Fischer Theorem) Let $\{a_n\}\in l^2.$ Then there is $f\in L^2(T)$ so
that $a_n(f)=a_n$ for all $n$ and $||f||_2^2=\sum |a_n|^2.$
\end{thm}

{\bf Proof:}\quad Define $f_N=\sum _{-N}^Na_n\chi^n.$ Then $f_N$ is Cauchy
sequence in $L^2(T).$ Let $f$ be the limit of $f_N$ in $L^2(T).$

We verify that $a_n(f)=a_n$ for all $n.$ Fix $n$ and let $N\geq n.$ Note
that $a_n=\int f_N\chi^{-n}.$ Then 
\[
|a_n(f)-a_n|=|\int f\chi^{-n}-\int f_N\chi^{-n}| \leq ||f-f_N||_2\rightarrow
0,\quad \mbox{as $N\rightarrow\infty$}. 
\]

To prove $||f||_2^2=\sum |a_n|^2$ we note that $f_N\rightarrow f$ in $L^2(T)$
and so $||f_N||_2\rightarrow ||f||_2.$ Since $||f_N||_2^2=\sum_{-N}^N|a_n|^2
\rightarrow \sum |a_n|^2,$ $||f||_2^2=\sum |a_n|^2.$

\begin{thm}
For every $f\in L^2(T),$ the Parseval relation 
\[
||f||_2^2=\sum |a_n(f)|^2 
\]
holds. Equivalently, for every $f\in L^2(T)$, $f_N\rightarrow f$ in $L^2(T)$
as $N\rightarrow \infty.$
\end{thm}

{\bf Proof:}\quad Clearly, the Parseval relation holds for all trigonometric
polynomials (because the Fourier series of any trigonometric polynomial is
itself). Therefore, $\mathcal{F}$ (Fourier transform that carries $f\in L^{2}$
to $\{a_{n}(f)\}\in l^{2}$) is an isometry from trigonometric polynomials in
the norm of $L^{2}(T)$ into $l^{2},$ and its range consists of all sequence $
\{a_{n}\}$ such that $a_{n}=0$ from some $n$ on. 

Note that the range is a
dense subset of $l^{2}$. If we prove that the family of trigonometric
polynomials is dense in $L^{2}(T),$ then 
$\mathcal{F}$ has a unique continuous extension, also denoted as $\mathcal{F}$,
to a linear isometry of all of $L^{2}(T)$ onto $l^{2}.$ This extension $
\mathcal{F}$ must be the Fourier transform. (To show this, one needs to prove
that for any $f\in L^{2}(T)$, $\mathcal{F}$ $(f)=\{a_{n}(f)\}.$) Therefore,
the Fourier transform is an isometry from $L^{2}(T)$ onto $l^{2}.$

\begin{thm}
$\{\chi^n\}$ is complete. More precisely, if all the Fourier coefficients of 
$f\in L^1$ are zero, then $f=0$ a.e.
\end{thm}

{\bf Proof:} See ([1])

\begin{thm}
(Mercer's Theorem) For any $f\in L^1(T)$, $a_n(f)\rightarrow 0$ as $
n\rightarrow \pm\infty.$
\end{thm}
{\bf Proof:}\quad If $f\in L^{2}(T)$ then the theorem follows from the
Bessel's inequality. If $f\in L^{1}(T)$, choose $f_{k}\in L^{2}(T)$ so that $
f_{k}\rightarrow f$ in $L^{1}(T).$ 
Then for each $n$, $|a_{n}(f_{k})-a_{n}(f)|=|a_{n}(f_{k}-f)|\leq
||f_{k}-f||_{1}\rightarrow 0$ as $k\rightarrow \infty .$ This shows that $
\{a_{n}(f_{k})\}$ converges to $\{a_{n}(f)\}$ as $k\rightarrow \infty $
uniformly in $n.$ Now, write $|a_{n}(f)|\leq
|a_{n}(f)-a_{n}(f_{k})|+|a_{n}(f_{k})|$ and for any $\epsilon >0,$ choose $k$
so that $|a_{n}(f)-a_{n}(f_{k})|<\epsilon $ for all $n.$ By fixing this $k$
and letting $n$ large enough, we get $|a_{n}(f_{k})|<\epsilon .$

\begin{thm}
If $f\in L^1(T)$ and ${\displaystyle \frac{f(e^{ix})}{x}}$ is integrable on $
(-\pi,\pi),$ then ${\displaystyle \sum_{n=-M}^{N} a_n(f)\rightarrow 0}$ as $
M,N\rightarrow \infty$ (independently).
\end{thm}

{\bf Proof:}\quad By the hypothesis, 
\[
g(e^{ix})=\frac{f(e^{2ix})}{sinx}\in L^{1}(T). 
\]
(Note: the behavior of $g$ near $\pm \pi $ is analogous to that of ${
\displaystyle\frac{f(e^{ix})}{x}}$ near $0$). Rewriting 
\[
f(e^{2ix})=\frac{(e^{ix}-e^{-ix})g(e^{ix})}{2i} 
\]
and integrating against $\chi ^{-2n}d\sigma $, we get 
\[
2ia_{n}(f)=a_{2n-1}(g)-a_{2n+1}(g),\quad \forall n. 
\]
Hence, (telescoping sum), as $M,N\rightarrow \infty ,$ 
\[
2i\sum_{-M}^{N}a_{n}(f)=a_{-2M-1}-a_{2N+1}\rightarrow 0. 
\]
(It is worth noting that the gist of the proof is considering $f(e^{2ix})$
and ending up with a telescoping sum.)

\begin{cor}
If $f\in L^{1}$ and $f$ satisfies Lipschitz condition at $e^{it},$ then the
Fourier series of $f$ converges to $f$ at that point. That is, $\sum
a_{n}(f)e^{int}\rightarrow f(e^{it})$.
\end{cor}

{\bf Proof:}\quad Without loss of generality, we may assume that $t=0$ and $
f(1)=0$ and show that $\sum a_{n}(f)\rightarrow 0.$

Assume that $f$ satisfies the Lipschitz condition at $e^{it}$, that is,
there is a neighborhood of $t$ so that for any $x$ in that neighborhood, $
|f(e^{ix})-f(e^{it})|\leq K|x-t|^{\alpha }$ for some $0<\alpha \leq 1.$ In
our case of $t=0$ and $f(1)=0$, this means that $|f(e^{ix})|\leq
K|x|^{\alpha },$ for $x$ close to $0.$ Therefore, ${\displaystyle\frac{
f(e^{ix})}{x}}$ is integrable on $(-\pi ,\pi ).$ Now the corollary follows
from the above theorem.

\begin{thm}
(Principle of localization) If $f, g\in L^1(T)$ and $f=g$ on some interval,
then at each interior point of the interval, their Fourier series are
equi-convergent.
\end{thm}

{\bf Proof:}\quad Note that 
\[
\sum_{-M}^N a_n(f)\chi^n-\sum_{-M}^N a_n(g)\chi^n =\sum_{-M}^N
a_n(f-g)\chi^n. 
\]
Since $f-g=0$ on some interval, $f-g$ satisfies Lipschitz condition at each
interior point of that interval. Therefore, 
\[
\sum_{-M}^N a_n(f-g)\chi^n\rightarrow 0. 
\]
This completes the proof.

\begin{thm}
Suppose $f(e^{ix})\in L^{1}(T)$ and ${\displaystyle\frac{f(e^{ix})+f(e^{-ix})
}{x}}$ is integrable on $(-\pi ,\pi )$. Show that 
\[
\sum_{-N}^{N}a_{n}(f)\rightarrow 0\,\mbox{as $N\rightarrow \infty$}. 
\]
\end{thm}

{\bf Proof:}\quad Let $g(e^{ix})$ be such that 
\[
f(e^{2ix})-f(e^{-2ix})=\frac{1}{2i}(e^{ix}-e^{-ix})g(e^{ix}). 
\]
Note that $g$ is integrable on $(-\pi ,\pi )$ by the hypothesis. Integrating
against $\chi ^{-2n}d\sigma $ we have 
\[
2i(a_{n}(f)-a_{-n}(f))=a_{2n-1}(g)-a_{2n+1}(g). 
\]
Adding up these equalities for $n=0,\pm 1,\cdots ,\pm N$ we have 
\[
2i\sum_{-N}^{N}(a_{n}(f)+a_{-n}(f))=4i
\sum_{-N}^{N}a_{n}(f)=a_{-2N-1}(g)-a_{2N+1}(g)\rightarrow 0,\quad \mbox{as}
\quad N\rightarrow \infty . 
\]

It is worth noting that, under the hypothesis of the theorem, it is not
necessarily true that $\sum_{-M}^{N}a_{n}(f)\rightarrow 0$ as $
N,M\rightarrow \infty $ independently. For example, let $f(e^{ix})=-1$ on $
(-\pi ,0)$ and $=1$ on $(0,\pi ).$ Then $a_{n}(f)=\frac{-i}{n\pi }$ for $n$
odd, $=0$ for $n$ even and 
\[
f(e^{ix})\sim \sum \frac{1}{\pi i}(\frac{1-(-1)^{n}}{n}). 
\]
Clearly, $\sum_{-N}^{N}a_{n}(f)=0$ for all $N.$ However, if $M=2N$ and $N=2k$
then, 
\[
\sum_{-N}^{M}a_{n}(f)=\frac{2}{\pi }\sum_{k}^{2k-1}\frac{1}{2n+1}>\frac{2}{
\pi }(k\frac{1}{4k-1})\not\rightarrow 0 
\]
as $k\rightarrow \infty .$ Note that $f$ (with the definition $f(1)=0$ is
the value of the midpoint in the gap at $x=0$) satisfies the hypothesis in
Theorem 1.10 but not the condition in Theorem 1.8.

\begin{cor}
Assume that $f\in L^{1}$ and $f$ satisfies 
symmetric Lipschitz condition at $e^{it},$ that is, there
is a neighborhood of $t$ so that for any $x$ in that neighborhood, 
\[
|f(e^{ix})+f(e^{-ix})-2f(e^{it})|\leq K|x-t|^{\alpha } 
\]
for some $0<\alpha \leq 1.$ Then the Fourier series converges to $f$ at that
point. That is, $\sum_{N}^{N}a_{n}(f)e^{int}\rightarrow f(e^{it})$ as $
N\rightarrow \infty $.
\end{cor}
{\bf Proof:} See ([2])

\begin{thm}
For $0<r< 1,$ ${\displaystyle f(e^{ix})=\frac{1}{1-re^{ix}}=\sum r^ne^{inx}}$
so that $a_n(f)=r^n.$
\end{thm}

\vspace{0.1in}

\section{Convolution}

In this section, for any function $f$ defined on $T,$ we write
the value of $f$ at $e^{ix}$ as $f(x).$ That is, we always read $f(x)$ as
the value of $f$ at $e^{ix}$.

\vspace{0.2in}

{\bf  Convolution on the group ${\bf T}$} 

\vspace{0.1in}

\begin{df}
For $f, g\in L^1(T),$ 
\[
f*g(x)=\int f(t)g(x-t)d\sigma(t)= \int f(x-t)g(t)d\sigma(t). 
\]
\end{df}

\begin{thm}[Continuity in $L^p(T)$]
If $f\in L^p(T), 1\leq p<\infty,$ then 
\[
\lim_{|h|\rightarrow 0}||f(x+h)-f(x)||_p=0. 
\]
\end{thm}

{\bf Proof:}\quad Let 
\[
C_{p}=\{f\in L^{p}(T):||f(x+h)-f(x)||_{p}\rightarrow 0,\quad \mbox{as}\quad
h\rightarrow 0\}. 
\]

Claim:

\begin{enumerate}
\item A finite linear combination of functions in $C_p$ is in $C_p.$ 

Proof: Let $f$ and $g$ be in $C_p$ and $a,b$ are two numbers. Then $
||(af(x+h)+bg(x+h))-(af(x)+bg(x))||_p \leq
|a|||f(x+h)-f(x)||_p+|b|||g(x+h)-g(x)||_p\rightarrow 0$ as $h\rightarrow 0.$
\vspace{0.1in}

\item If $f_k\in C_p$ and $f_k\rightarrow f$ in $L^p,$ then $f\in C_p.$

Proof: Note that $||f(x+h)-f(x)||_p\leq
||f(x+h)-f_k(x+h)||_p+||f_k(x+h)-f_k(x)||_p+
||f_k(x)-f(x)||_p=||f_k(x+h)-f_k(x)||_p+2||f(x)-f_k(x)||_p.$ Since $f_k\in
C_p,$ we have $\limsup_{|h|\rightarrow 0}||f(x+h)-f(x)||_p \leq
2||f_k(x)-f(x)||_p,$ and this goes to zero by letting $k\rightarrow \infty.$
\end{enumerate}

Clearly, the characteristic function of an interval belongs to $C_{p}.$
Since the linear combinations of characteristic functions (step functions)
of intervals are dense in $L^{p}(T),$ by using the method of successively
approximating more and more general functions, it follows that $L^{p}(T)$ is
contained in $C_{p}.$

\begin{thm}
\begin{enumerate}
\item If $f\in L^p(T), 1\leq p\leq \infty,$ and $g\in L^{p^{\prime}}(T),$
where $\frac{1}{p} +\frac{1}{p^{\prime}}=1,$ then $f*g(x)$ exists everywhere
and is continuous with 
\[||f*g||_C\leq ||f||_p\cdot ||g||_{p^{\prime}}. \]

\vspace{0.1in}

\item If $f\in L^p, 1\leq p<\infty,$  and $g\in
L^1(T), $ then $f*g(x)$ exists a.e. as an absolutely convergent integral, $f*g\in L^p(T),$ and $||f*g||_\leq ||f||_p\cdot
||g||_1.$ 
\vspace{0.1in}
\item If $f\in C(T)$ and $g\in
L^1(T), $ then $f*g(x)$ exists a.e. as an absolutely convergent integral, $f*g\in C(T)$, and  $||f*g||_C\leq ||f||_C\cdot ||g||_1.$

\end{enumerate}
\end{thm}

{\bf Proof:}\quad

\begin{enumerate}
\item $f\in L^p, 1\leq p<\infty$ and $g\in L^{p^{\prime}}.$

Proof: Using the H\"{o}lder inequality and continuity in $L^{p},$ we see that 
\begin{eqnarray*}
|f\ast g(x+h)-f\ast g(x)| &\leq &\int |f(x+h-t)-f(x-t)||g(t)|d\sigma (t) \\
&\leq &||f(x+h-\cdot )-f(x-\cdot )||_{p}\cdot ||g||_{p^{\prime }} \\
&\leq &||f(\cdot +h)-f(\cdot )||_{p}\cdot ||g||_{p^{\prime }}
\end{eqnarray*}
approaches $0$ as $h\rightarrow 0.$ The last inequality above holds because
the measure $\sigma $ is translation invariant. The norm estimate follows
from the H\"{o}lder inequality again. If $p=\infty ,$ the roles of $f$ and $g $ may be interchanged.

\vspace{0.1in}

\item $f\in L^p(T), 1\leq p<\infty,$ and $g\in L^1$.

Proof: Since for almost all $u$ 

\[\int |f(x-u)|^p|g(u)|d\sigma(x)=|g(u)|\cdot ||f||_p^p 
\]
which belongs to $L^1(T),$ it follows that 
\[
\int\int |f(x-u)|^p|g(u)|d\sigma(x)d\sigma(u) =||f||_p^p||g||_1
\]
exists as a finite number. Therefore by Fubini's theorem 
\[
\int\biggl(\int |f(x-u)|^p|g(u)|d\sigma(u)\biggr)d\sigma(x)  
\]
exists as well and is equal to $||f||_p^p||g||_1.$ This implies that 
\[
\int |f(x-u)|^p|g(u)|d\sigma(u) 
\]
exists for almost every $x\in R$ and belongs to $L^1(T).$ This proves the
assertion for $p=1.$

For $1<p<\infty,$ H\"{o}lder's inequality delivers 
\[
|(f\ast g)(x)|\leq \biggl(\int |f(x-u)|^{p}|g(u)|d\sigma (u)\biggr)^{1/p}
\biggl(\int |g(u)|d\sigma (u)\biggr)^{1/p^{\prime }}. 
\]
This shows that $f\ast g(x)$ exists a.e. as an absolutely
convergent integral. Moreover, 
\begin{eqnarray*}
||f\ast g||_{p} &\leq &||g||_{1}^{1/p^{\prime }}\cdot \biggr(\int \int
|f(x-u)|^{p}|g(u)|d\sigma (u)d\sigma (x)\biggl)^{1/p} \\
&=&||g||_{1}^{1/p^{\prime }}||f||_{p}||g||_{1}^{1/p}=||f||_{p}||g||_{1}
\end{eqnarray*}

(3). Finally, if $f\in C(T),$ it follows as in the proof of (1) that $f*g(x)$
exists for every $x,$ belongs to $C(T)$ and satisfies $||f*g||_C \leq
||f||_C||g||_1.$
\end{enumerate}

\begin{df}
A (complex) linear algebra is a (complex) linear space $\mathcal{A}$ in which a product is defined such that, for all $x,y,z\in \mathcal{A}, a,
b\in C$ $x(yz)=(xy)z$ (associative), $(ax)y=x(ay)=a(xy)$ (commutative with
scalar), $x(ay+bz)=a(xy)+b(xz)$ (distributive over addition, or say linear).

If a linear algebra $\mathcal{A}$ is equipped with a norm under which it is a Banach space, $\mathcal{A}$ is a Banach algebra if $||xy||\leq
||x||\cdot ||y||$ for all $x,y\in \mathcal{A}.$ If furthermore $\mathcal{A}$
has an identity for multiplication, $e=ex=xe$ for all $x\in \mathcal{A},$ then 
$||e||=1.$

\end{df}

\begin{thm}
$L^1(T)$ is a commutative Banach algebra under convolution. This algebra has
no identity.
\end{thm}

{\bf Proof:}\quad We show that this algebra has no identity. If there were
such function $\phi\in L^1(T)$ it would mean that $(\phi*f)(x)=f(x)$ for all 
$f\in L^1(T)$ and almost every $x\in R.$ Looking at $\phi*f(x)=\int
\phi(x-t)f(t)d\sigma(t),$ since the value of $f$ at $x$ (i.e. $e^{ix}$) does
not depend on the values $f$ takes at other points, we see that $\phi(x-t)=0$
for almost all $t\neq x.$ Otherwise, changing $f$ at $t$ for $t$ in a set of
positive measure would change the value of $f(x).$ So $\phi=0$ a.e., but $
\phi*f$ is not identically $0$ if $f$ is not identically zero.

{\bf Convolution on the line group ${\bf R}$}

\begin{df}
Let $f, g$ be two functions defined and measurable on $R$. The expression 
\[
f*g(x)=\int f(t)g(x-t)d\sigma(t) =\int f(x-t)g(t)d\sigma(t) 
\]
is called the convolution of $f$ and $g$.
\end{df}

\begin{thm}
\begin{enumerate}
\item Let $f\in L^p, 1\leq p\leq \infty,$ and $g\in L^{p^{\prime}}.$ Then $
f*g(x)$ exists everywhere, belongs to $C$, and $||f*g||_C \leq
||f||_p||g||_{p^{\prime}}.$ Moreover, if $1<p<\infty$, then $f*g\in C_0$,
i.e. $f*g\in C$ and $\lim_{|x|\rightarrow \infty} f*g(x)=0.$ The same is
true for $p=1,$ if, in addition, $g\in C_0.$
\vspace{0.1in}
\item If $f\in L^p, 1\leq p<\infty,$ and $g\in L^1,$
then $f*g(x)$ exists a.e. as an absolutely convergent integral, $f*g\in L^p,$ and $||f*g||_p\leq ||f||_p\cdot ||g||_1.$
\vspace{0.1in}
\item If $f\in C$ and $g\in L^1,$
then $f*g(x)$ exists a.e. as an absolutely convergent integral,  $f*g\in C$, and  $||f*g||_C\leq ||f||_C\cdot ||g||_1.$
\end{enumerate}
\end{thm}
{\bf Proof:} (See [2])

\begin{df}
An approximate identity on $T$ is a sequence of functions $e_n$
with these properties:

\begin{enumerate}
\item each $e_n$ is nonnegative;

\item ${\displaystyle\int e_n(t)d\sigma(t)=1};$

\item for every $0<\delta<\pi,$ 
\[
\lim_{n\rightarrow\infty}\int_{|t|>\delta}e_n(t)d\sigma(t)=0. 
\]
\end{enumerate}
\end{df}

\begin{thm}
Let $X$ and $Y$ be Banach spaces, and $T_n$ a sequence of linear operators
from $X$ to $Y$ with $||T_n||\leq K$ for all $n,$ and let $D$ be a dense
subset of $X.$ Suppose for each $x\in D,$ $T_nx$ converges. Then $T_nx$
converges for all $x\in X$ and if we define $T: X\rightarrow Y$ by $Tx=lim
T_nx$ for all $x\in X,$ then $T$ is a linear operator with $||T||\leq K.$
\end{thm}

{\bf Proof:}\quad Let $x\in X$ and $x_k\in D$ with $x_k\rightarrow x$ in $X.$
Then $||T_nx-T_mx||\leq ||T_nx-T_nx_k||+||T_nx_k-T_mx_k||+||T_mx_k-T_mx||
\leq 2K||x-x_k||+||T_nx_k-T_mx_k||.$ Given $\epsilon>0$, find $k$ so that $
2K||x-x_k||<\epsilon$ and for this fixed $k$, find $N$ so that if $n,m>N,$
then $||T_nx_k-T_mx_k||<\epsilon.$ So $T_nx$ is a Cauchy sequence in $Y$.
Since $Y$ is Banach, $T_nx$ converges.

Let $T$ be defined as in the theorem. Then $T$ is linear operator and for
any $x$ with $||x||\leq 1,$ $||Tx||\leq K.$ Therefore, $||T||=\sup_{||x||
\leq 1}||Tx||\leq K.$

\begin{thm}[Fej\'{e}r's]
Let $f\in L^p(T)$ with $1\leq p<\infty.$ Then for any approximate identity $
e_n,$ $e_n*f\rightarrow f$ in $L^p(T).$ If $f\in C(T),$ then $
e_n*f\rightarrow f$ uniformly.
\end{thm}

{\bf Proof:}\quad If $f\in C(T),$ then 
\[
f(x)-e_n*f(x)=\int (f(x)-f(x-t)) e_n(t)d\sigma(t). 
\]
Note that $f$ is uniformly continuous on $T$. Given $\epsilon>0$, there is $
\delta>0$ such that $|f(x)-f(x-t)|<\epsilon$ for $|t|\leq \delta.$ Denote by 
$I$ the part of the integral over $|t|<\delta$, and by $J$ the integral over
the complementary interval. Clearly, 
\[
|I|\leq \epsilon \int_{|t|\leq \delta}e_n(t)d\sigma(t) \leq \epsilon 
\]
and 
\[
|J|\leq 2M\int_{|t|>\delta}e_n(t)d\sigma(t)\rightarrow 0 
\]
as $n\rightarrow \infty.$

Now let $f\in L^p$ and $1\leq p<\infty.$ For continuous $f,$ the uniform
convergence of $e_n*f$ to $f$ implies convergence in $L^p.$ Let $
T_n(f)=e_n*f.$ Then $T_n$ a linear operator from $L^p(T)\rightarrow L^p(T)$
such that $||T_n(f)||_p\leq ||e_n||_1||f||_p$, i.e., $||T_n||\leq 1.$ Note
that $T_n(f)$ converges in $L^p(T)$ for every $f\in C(T)$ and $C(T)$ is
dense in $L^p(T).$ Therefore, Theorem 1.4 asserts that $T_n(f)$ converges in 
$L^p(T)$ for every $f\in L^p(T)$ and if we define $T(f)=L^p-\lim T_n(f)$
then $T$ is a linear operator on $L^p$ with bound $\leq 1.$ We prove that $T$
is an identity on $L^p$. In fact, $T(f)=f$ for all $f\in C(T)$ and $C(T)$ is
dense in $L^p(T).$ Let $f\in L^p$ and let $f_k\in C(T)$ with $f_k\rightarrow
f$ in $L^p.$ Then $T(f)=\lim T(f_k)=\lim f_k=f$ for all $f$, that is, $
L^p-\lim e_n*f=f$ for all $f\in L^p.$

\vspace{0.1in}

\begin{thm}
Let $e_n$ be an approximate identity on $T$, and $g\in L^\infty(T).$ Then $
e_n*g\rightarrow g$ in the weak* topology of $L^\infty.$ That is, for every $
f\in L^1(T),$ 
\[
\lim_{n\rightarrow\infty}\int \biggl((e_n*g)(x)-g(x)\biggl) f(x)
d\sigma(x)\rightarrow 0. 
\]
\end{thm}

{\bf Proof:}\quad If $g\in L^\infty$ then 
\[
e_n*g(x)-g(x)=\int (g(x-u)-g(x))e_n(u)d\sigma(u). 
\]
By Fubini's theorem, for every $f\in L^1(T),$ and $0<\delta<\pi,$ 
\begin{eqnarray*}
&.&\int \biggl((e_n*g)(x)-g(x)\biggl) f(x) d\sigma(x) \\
&=&\int \biggl(\int (g(x-u)-g(x))e_n(u)d\sigma(u)\biggr) f(x) d\sigma(x) \\
&=&\int \int (g(x-u)f(x)-g(x) f(x)) d\sigma(x) e_n(u)d\sigma(u) \\
&=&\int \int (g(x)f(x+u)-g(x) f(x)) d\sigma(x) e_n(u)d\sigma(u) \\
&=&\int \int (f(x+u)-f(x))g(x)d\sigma(x) e_n(u)d\sigma(u).
\end{eqnarray*}
But since $f\in L^1(T)$ is continuous in mean, for each $\epsilon>0,$ there
is a $\delta>0$ such that $||f(\cdot+u)-f(\cdot)||_1<\epsilon$ for all $
|u|\leq \delta.$ Denote by $I$ the part of the integral over $|u|<\delta$,
and by $J$ the integral over the complementary interval. Then 
\[
|I|\leq sup_{|u|\leq \delta} |\int g(x)(f(x+u)-f(x)) d\sigma(x)| \leq
||g||_\infty \epsilon. 
\]
Furthermore, 
\[
|J|\leq 2||g||_{\infty}||f||_1\int_{\delta\leq |u|\leq \pi}|e_n(u)|
d\sigma(u). 
\]
It follows $J\rightarrow 0$ as $n\rightarrow\infty.$

\begin{thm}
Let $e_k$ be an approximate identity on $T$. Then for each $n,$ $
a_n(e_k)\rightarrow 1$ as $k\rightarrow \infty.$
\end{thm}

{\bf Proof:}\quad For each $n,$ $e_k*e^{in\cdot}(e^{ix})\rightarrow e^{inx}$
uniformly. In particular, let $x=0,$ $e_k*e^{in\cdot}(1)\rightarrow 1.$ That
is, $a_n(e_k)\rightarrow 1.$

\vspace{0.1in}

\begin{thm}
Assume that there is an approximate identity $e_k$ on $T$ consisting of
trigonometric polynomials. Then trigonometric polynomials are dense in $
L^p(T)$, $1\leq p<\infty.$
\end{thm}

{\bf Proof:}\quad Let $f\in L^p,$ $1\leq p<\infty.$ Then $f\in L^1(T)$ and $a_n(e_k*f)=a_n(e_k)a_n(f).$ Since $e_k$ is a trigonometric polynomial, for
each $k$, $a_n(e_k)=0$ so that $a_n(e_k*f)=0$ for large $n.$ Thus, $e_k*f$
is a trigonometric polynomial. Since $e_k*f\rightarrow f$ in $L^p(T),$ 
$1\leq p<\infty,$ trigonometric polynomials are dense in $L^p(T).$

\begin{thm}[Converse of H\"{o}lder's Inequality]
Let $1\leq p<\infty.$ If $\int |fg|d\sigma\leq k$ for every $g\in
L^{p^{\prime}}(T)$ with norm $1,$ then $f\in L^{p}(T)$ and $||f||_p\leq k.$
\end{thm}

{\bf Proof:}\quad Let $f_n=f,$ if $x\in\{x:|f(x)|\leq n\};$ $f_n=0$ if $x$
is in the complement. Then $f_n\in L^p(T)$ for any $p$ and $|f_n|\uparrow
|f| $ a.e. Let 
\[g=\frac{|f_n|^{p-1}}{||f_n||_p^{p-1}}. \]
Then $||g||_{p^{\prime}}=1.$ By assumption, we have 
\[\frac{1}{||f_n||_p^{p-1}}\int |f||f_n|^{p-1} d\sigma=\frac{1}{||f_n||_p^{p-1}} \int |f_n|^pd\sigma\leq k. \]
By B. Levi's theorem (Monotone Convergence Theorem), $||f||_p \leq k.$ The
proof can be extended to $L^p(R)$ by defining $f_n=f$ for $x\in \{x\in R:
|x|\leq n \quad\mbox{and}\quad |f(x)|\leq n\}.$

\section{Unicity Theorem, Parseval Relation}

\begin{thm}[Unicity Theorem]
If $f\in L^1(T)$ and $a_n(f)=0$ for all $n,$ then $f=0$ a.e.
\end{thm}

{\bf Proof:}\quad If $f$ has continuous derivative then $f$ satisfies the Lipschitz condition at every point so that $\sum a_n(f)e^{inx}$ converges to 
$f(x)$ everywhere. If $a_n(f)=0$ for all $n,$ then the sum is zero so that $
f=0$ everywhere.

Let $e_n$ be an approximate identity in $L^1(T)$ so that $e_n$'s are
continuously differentiable. Then for every $f\in L^1,$ $e_n*f$ is continuously
differentiable. Let $f$ be such that $a_n(f)=0$ for all $n.$ Then for every $
n,$ $a_k(e_n*f)=a_n(e_k)a_n(f)=0$ for all $k.$ By the first part of proof, $
e_n*f=0$ everywhere. Thus, $f,$ as a limit of $e_n*f$ in $L^1,$ is zero
almost everywhere.

\begin{cor}
The trigonometric polynomials are dense in $L^p(T), 1\leq p<\infty.$
\end{cor}

{\bf Proof:}\quad If the trigonometric polynomials are not dense in $L^p(T),$
then there exists a nonzero linear functional $l$ on $L^p(T)$ so that $l(h)=0 $ for all $h$ in the closure of the trigonometric polynomials. It follows that there
exists $g\in L^{p^{\prime}}, 1<p^{\prime}\leq \infty,$ $g\neq 0, a.e.$ so
that $l(h)=\int g\overline{h}d\sigma=0$ for all trigonometric polynomials $
h. $ This means $a_n(g)=0$ for all $n$ so that $g=0$ a.e., contrary to
hypothesis.

\begin{thm}[Parseval Relation]
The Parseval relation 
\[
||f||_2^2=\sum |a_n(f)|^2 
\]
holds for every $f\in L^2(T).$
\end{thm}

{\bf Proof:}\quad The Parseval relation holds for all trigonometric
polynomials. It is also valid for every function that is the limit of
trigonometric polynomials in $L^2(T).$ Since every function in $
L^2(T) $ is such a limit, the Parseval relation holds for every $f\in
L^2(T). $

\vspace{0.1in}

Let $C_0(X)$ denote the subset of $C(X)$ consisting of all $f\in C(X)$ such
that for every $\epsilon>0$, there is a compact subset $F$ of $X$ such that $
|f(x)|<\epsilon$ for all $x\in F^{\prime}\bigcap X.$ If $X$ is compact then $
C_0(X)=C(X).$

\begin{thm}[F. Riesz Representation Theorem]
Let $X$ be a locally compact Hausdorff space. Then the mapping $T(\mu)=l_\mu$
is a norm-preserving linear mapping of $M(X)$ onto $C_0(X)^*,$\quad where $
l_\mu(f)=\int_X fd\mu,$ $f\in C_0(X).$ Thus $M(X)$ is a Banach space, and $
M(X)$ and $C_0(X)^*$ are isomorphic as Banach spaces.
\end{thm}
{\bf Proof:} See ([3])

\begin{thm}[Weierstrass]
The trigonometric polynomials are dense in $C(T).$
\end{thm}

{\bf Proof:}\quad We want to show that any continuous linear functional $l$
on $C(T)$ that vanishes on all trigonometric polynomials is null. Note that such a linear functional can be realized by a
bounded, complex Borel measure $\mu$ on $T$ (realized as $[-\pi, \pi)$) in
the way that
\[l(h)=\int h(e^{-ix})d\mu(x),\qquad h\in C(T).\]
Thus, the proof of Weierstrass' theorem is converted to proving that if $\int
e^{-inx}d\mu(x)=0$ for all $n$ then $\mu$ is null. We prove this fact in the
following theorem, which is a stronger statement of Unicity theorem when we
identify a function with the measure $f(e^{ix})d\sigma(x).$

\begin{thm}[Unicity Theorem]
Let $\mu$ be Borel measure on $T.$ The Fourier-Stieltjes coefficients of $
\mu $ are defined to be 
\[
\hat{\mu}(n)=a_n(\mu)=\int e^{-inx}d\mu(x). 
\]
If $\hat{\mu}(n)=0$ for all $n,$ then $\mu$ is null.
\end{thm}

{\bf Proof:}\quad Note that if $f\in C(T)$ and $a_n(f)=0$ for all $n$ then $
f=0.$ This is the weakest Unicity theorem of all. We use this version and
the Riesz theorem to prove the strongest version as stated in the theorem.

For $g\in C(T),$ define the convolution 
\[
\mu *g(e^{ix})=\int g(e^{i(x-t)})d\mu(t). 
\]
We show that it is a continuous function.  As usual, we write $g(e^{ix})$
as $g(x)$ for simplicity. Note that 
\begin{eqnarray*}
|\mu *g(x+h)-\mu*g(x)|&=&|\int (g(x+h-t)-g(x-t))d\mu(t)| \\
&\leq& ||g(\cdot+h)-g(\cdot)||_C\int d|\mu(t)| \\
&=& ||g(\cdot+h)-g(\cdot)||_C||\mu||.
\end{eqnarray*}
Since $|\mu|(T)$ is finite and $g$ is uniformly continuous, the last
expression tends to zero as $h\rightarrow 0.$

The Fourier coefficients of $\mu*g$ are $a_n(\mu *g)=a_n(\mu)a_n(g).$ By
assumption, $a_n(\mu)=0$ for all $n.$ Thus, $a_n(\mu *g)=0$ for all $n$ and,
since $\mu*g\in C(T),$ $\mu*g(x)=0$ at every $x.$ In particular, $
\mu*g(1)=0. $ Since $g$ was an arbitrary function in $C(T),$ $(\mu* g)(1)$
is the null functional on $g\in C(T),$ and so by the Riesz theorem, $\mu$ is
null measure. Note that $\mu*g(1)=\int g(e^{-ix})d\mu(x).$ The theorem is
proved.

\begin{df}
The convolution of two Borel measures on $T$, say $\mu$ and $\nu$, is
defined to be the measure $\mu*\nu$ such that 
\[
(\mu*\nu)*h(1)=(\mu*(\nu*h))(1). 
\]
\end{df}

Clearly, by translating $h$ we see that $(\mu*\nu)*h(e^{ix})=(\mu*(
\nu*h))(e^{ix}),$ for $h\in C(T)$ and all $x.$ That is, for all $x,$ 
\[
\int h(e^{i(x-t)})d(\mu*\nu)(t)= =\int \int h(e^{i(x-s-t)})d\nu(t)d\mu(s). 
\]
Letting $x=0$ gives the equality in the above definition of convolution of
two measures.

\begin{thm}
$\mu*\nu\in M(T).$
\end{thm}

{\bf Proof:}\quad Define the functional $l(h)=(\mu*(\nu*h))(1)$ on $C(T).$
Clearly, it is linear. Observe that $||\nu*h||_\infty\leq ||\nu||
||h||_\infty$ for any measure $\nu\in M(T)$ and $h\in C(T).$ Applying twice,
we have $|l(h)|\leq ||\mu|| ||\nu|| ||h||_\infty.$ Therefore, $l$ is a
continuous linear functional on $C(T).$ By the Riesz theorem, there is $
\gamma\in M(T)$ so that $l(h)=\int h(e^{-it})d\gamma(t)=\gamma*h(1)$ for all 
$h\in C(T).$ Then $\mu*\nu$ is defined to be the measure $\gamma\in M(T).$
Moreover, $||\mu*\nu||\leq ||\mu|| ||\nu||.$

\begin{thm}
$M(T)$ is a commutative Banach algebra under the convolution.
\end{thm}

{\bf Proof:}\quad Observe that, for $\chi$ defined as $\chi(e^{ix})=e^{ix},$ 
$(\mu* \chi^n)(e^{ix})=\int e^{in(x-t)}d\mu(t)=a_n(\mu)\chi^n(e^{ix}).$
Hence, 
\begin{eqnarray*}
a_n(\mu*\nu)\cdot\chi^n(1)&=& ((\mu*\nu)*\chi^n)(1)= (\mu*(\nu*\chi^n))(1) \\
&=&a_n(\nu) (\mu*\chi^n)(1)=a_n(\nu)a_n(\mu)\cdot \chi^n(1).
\end{eqnarray*}
It follows from the Unicity theorem that $\mu*\nu=\nu*\mu*,$ and $
(\mu*\nu)*\gamma=\mu*(\nu*\gamma).$

\begin{prob}
The convolution of a measure with an integrable function is the same as that
induced by considering $L^1(T)$ as a subalgebra of $M(T).$ In other words, $
\mu*(f d\sigma)=(\mu*f)d\sigma,$ where $f\in L^1(T)$ and is identified with
the measure $f d\sigma,$ while $\mu*f\in L^1(T).$
\end{prob}

{\bf Solution:}\quad We have shown that if $f\in C(T),$ then $\mu*f\in C(T).$
Let $f\in L^p,$ $1\leq p<\infty.$ It follows by Fubini's theorem that 
\[
\int (\int |f(e^{i(x-t)})|^pd|\mu(t)|)d\sigma(x)=||f||_p||\mu||. 
\]
This implies that $\int |f(e^{i(x-t)}|^pd|\mu(t)|$ exists for almost all $x$
and belongs to $L^1(T).$ This proves that if $f\in L^1(T),$ then so is $
f*\mu. $ If $1<p<\infty,$ by H\"{o}lder's inequality 
\[
|f*\mu(e^{ix})|\leq \int |f(e^{i(x-t)})d|\mu(t)| \leq (\int
|f(e^{i(x-t)})|^pd|\mu(t)|)^{1/p} (\int d|\mu(t)|)^{1/p^{\prime}}. 
\]
Since $\mu$ is bounded, $(f*\mu)(e^{ix})$ exists a.e. Further, one can show
that $f*\mu\in L^p(T)$ as the proof of the assertion that if $f\in L^p$ and $
g\in L^1,$ then $f*g\in L^p.$

We verify that both $\mu*(f d\sigma)$ and $(\mu*f)d\sigma$ are in $M(T)$
and, as linear functionals acting on $h\in C(T)$, 
\[
\int h(e^{-it})d(\mu*f d\sigma)(t) =\int h(e^{-it})(\mu*f)(e^{it})
d\sigma(t). 
\]
Thus $\mu*(f d\sigma)=(\mu*f)d\sigma.$

Note that the right hand side equals, using Fubini's theorem, 
\begin{eqnarray*}
\int h(e^{-it})\int f(e^{i(t-u)})d\mu(u) d\sigma(t)&=& \int \biggl(\int
h(e^{-it})f(e^{i(t-u)})d\sigma(t)\biggr)d\mu(u) \\
&=&\int (f*h)(e^{-iu})d\mu(u).
\end{eqnarray*}
The left hand side, by the definition of the convolution of two measures, equals 
\begin{eqnarray*}
\int h(e^{-it})\int f(e^{i(t-u)})d\mu(u)d\sigma(t)&=& \int \int
h(e^{-i(t+u)})f(e^{it})d\sigma(t)d\mu(u) \\
&=&\int (f*h)(e^{-iu})d\mu(u).
\end{eqnarray*}

\vspace{0.1in} The following theorem shows that the discrete part of a
measure $\mu$ can be extracted from its Fourier-Stieltjes series:

\begin{thm}[Wiener]
If $\mu\in M(T)$ has (its discrete part) point masses $a_n,$ then 
\[
\lim_{N\rightarrow \infty}\frac{1}{2N+1}\sum_{-N}^N|\hat{\mu}(n)|^2 =\sum
|a_n|^2. 
\]
Consequently, $\mu\in M(T)$ is continuous if and only if 
\[
\lim_{N\rightarrow \infty}\frac{1}{2N+1}\sum_{-N}^N|\hat{\mu}(n)|^2=0. 
\]
\end{thm}

{\bf Proof:}\quad Let $\tilde{\mu}(E)=\overline{\mu}(-E)$ for all Borel sets 
$E$ of $T.$ Then 
\begin{eqnarray*}
\hat{\tilde{\mu}}(n)&=&\int e^{-inx}d\tilde{\mu}(x)= \int e^{inx}d\tilde{\mu}(-x) \\
&=&\int e^{inx}d\overline{\mu}(x)= \overline{\int e^{-inx}d\mu(x)}=\overline{\hat{\mu}}
\end{eqnarray*}
Hence $a_n(\mu*\tilde{\mu})=a_n(\mu)a_n(\tilde{\mu})=|a_n(\mu)|^2$ so that 
\[
\frac{1}{2N+1}\sum_{-N}^N|\hat{\mu}(n)|^2=\int \frac{1}{2N+1}
\sum_{-N}^Ne^{-inx}d(\mu*\tilde{\mu})(x). 
\]
As $N\rightarrow \infty,$ the integrand tends boundedly to $1$ at $x=0$, and
to $0$ elsewhere on $(-\pi, \pi).$ Hence the limit is the mass of $\mu*
\tilde{\mu}$ at $0.$ 

To prove $\mu*\tilde{\mu}(\{0\})=\sum |a_n|^2,$ we define 
\[
h_\epsilon(e^{it})=\left\{
\begin{array}{ll}
1 & \mbox{if $|t|<\epsilon$,} \\ 
0 & \mbox{otherwise}.
\end{array}
\right. 
\]
Then, using the Bounded Convergence Theorem we can write 
\[
(\mu*\tilde{\mu})(\{0\})=\lim_{\epsilon\uparrow 0} \int h_\epsilon
(e^{-it})d\mu*\tilde{\mu}(t), 
\]
which can be viewed as the value of $(\mu*\tilde{\mu})*h_\epsilon (e^{ix})$
at $x=0.$ By the definition of the convolution of two measures, we can
rewrite 
\begin{eqnarray*}
(\mu*\tilde{\mu})*h_\epsilon (1)&=& \mu*(\tilde{\mu}*h_\epsilon)(1) \\
&=&\int \int h_\epsilon(e^{i(-s-t)})d\tilde{\mu}(t)d\mu(s) \\
&=&\int \int h_\epsilon(e^{i(-s+t)})d\tilde{\mu}(-t)d\mu(s) \\
&=&\int \int h_\epsilon(e^{i(t-s)})d\overline{\mu}(t)d\mu(s).
\end{eqnarray*}

We calculate the last convolution. For each $e^{is}\in T,$ 
\[
\lim_{\epsilon\rightarrow 0} \int h_\epsilon(e^{i(t-s)})d\overline{\mu}(t)
=\lim_{\epsilon\rightarrow 0} \int_Ad\overline{\mu}(t)=\lim_{\epsilon
\rightarrow 0}\overline{\mu}(A) =\overline{\mu}(\{e^{is}\}), 
\]
where $A=\{e^{it}\in T: s-\epsilon<t<s+\epsilon\}$ shrinks to $\{e^{is}\}$
as $\epsilon\rightarrow 0.$ Also note that $\mu\in M(T)$ so that $|\mu|(T)$
is finite. It follows that 
\[
|\int h_\epsilon(e^{i(t-s)})d\overline{\mu}(t)| \leq \int d|\overline{\mu}
(t)|=|\mu|(T)<\infty. 
\]
Therefore, by the Lebesgue Bounded Convergence Theorem, 
\begin{eqnarray*}
\int \int h_\epsilon(e^{i(t-s)})d\overline{\mu}(t)d\mu(s) &\rightarrow& \int 
\overline{\mu}(\{s\})d{\mu}(s) \\
&=&\sum \overline{\mu}(\{e^{i\theta_k}\}) \mu(\{e^{i\theta_k}\})=sum_k
|a_k|^2,
\end{eqnarray*}
as $\epsilon\rightarrow 0.$

\begin{prob}
If $\mu$ is a measure and $e_n$ is an approximate identity on $T$, then $
\mu*e_n$ tends to $\mu$ in the weak* topology of $M(T)$ as the dual of $
C(T). $ That is, 
\[
\lim_{n\rightarrow \infty}\int h(e^{-ix})\int
e_n(e^{i(x-t)})d\mu(t)d\sigma(x)\rightarrow \int h(e^{-it})d\mu(t) 
\]
for all $h\in C(T).$
\end{prob}

{\bf Proof:}\quad By Fubinis' theorem, 
\[
\int h(e^{-ix})\int e_n(e^{i(x-t)})d\mu(t)d\sigma(x) =\int \int
h(e^{-ix})e_n(e^{i(x-t)})d\sigma(x) d\mu(t). 
\]
Note that 
\[
\int h(e^{-ix})e_n(e^{i(x-t)})d\sigma(x) =e_n*h(e^{-it})\rightarrow
h(e^{-it}) 
\]
uniformly on $T.$ Taking limit under the outer integral gives rise to the desired limit.

\section{The Classical Kernels}

{\bf The Dirichlet kernel}

For each $n=0,1,2,\cdots,$ the partial sum $S_n(f)$ of Fourier series of a
function $f$ can be written as
\begin{eqnarray*}
S_n(f)(e^{ix})=\sum_{-n}^na_k(f)e^{ikx}&=& \sum_{-n}^n\int
f(e^{it})e^{ik(x-t)}d\sigma(t) \\
&=&\int f(e^{i(x-t)})\sum_{-n}^ne^{ikt}d\sigma(t) = (f*D_n)(e^{ix}),
\end{eqnarray*}
where 
\[D_n(e^{ix})=\sum_{-n}^ne^{ikx} =\frac{\sin (n+\frac{1}{2})x}{\sin \frac{x}{2}}, n=0,1,2,\cdots. 
\]

Since $D_n(e^{-ix})=D_n(e^{ix})$ (even), the Fourier series $S_n(f)$ of $f$
converges at $x=0$ as $n\rightarrow \infty$ if and only if 
\begin{eqnarray*}
\lim_{n\rightarrow \infty}(f*D_n)(1)&=&\lim_{n\rightarrow \infty} \int
f(e^{-it})D_n(e^{it})d\sigma(t) \\
&=&\lim_{n\rightarrow \infty} \int f(e^{it})D_n(e^{it})d\sigma(t)
\end{eqnarray*}
exists. If we take $f(e^{ix})=-1$ on $(-\pi,0)$ and $=1$ on $(0,\pi),$ then 
\[f\sim \sum_{-\infty}^\infty \frac{1}{\pi i}(\frac{1-(-1)^n}{n}
)e^{inx}.\]
 We have $(f*D_n)(1)=0$ for all $n$, but $\sum_{-M}^Na_n(f)$ does
not converge. Hence, the ''convergence'' of a Fourier series in complex form should be the convergence of symmetric sums.

\begin{thm}
If we define $L_n=||D_n||_1$ as the Lebesgue constant, then $L_n= \frac{4}{
\pi^2}\ln n+O(1).$ Thus $D_n$ do not form an approximate identity.
\end{thm}

{\bf Proof:}\quad Note 
\[
||D_n||_1= \int_{-\pi}^{\pi} |D_n|d\sigma=\frac{1}{\pi}\int_0^\pi |D_n|dx. 
\]
Since 
\[
\int_{0}^{\pi} |\sin(n+\frac{1}{2})x| |\frac{1}{\sin \frac{x}{2}}-\frac{1}{
\frac{x}{2}}|dx 
\]
is majorized by an absolute constant, we estimate 
\[
2\int_0^\pi |\frac{\sin(n+\frac{1}{2})x}{x}|dx, 
\]
which can be written as, if let $(n+\frac{1}{2})x=u,$ 
\[
2\int_0^{(n+\frac{1}{2})\pi} \frac{|\sin u|}{u}du. 
\]
We may disregard the parts of this integral over $(0,\pi)$ and $(n\pi, (n+
\frac{1}{2}\pi),$ since the integrand is bounded. In view of the periodicity
of $\sin u,$ what remains can be written as 
\[
2\int_\pi^{n\pi} \frac{|\sin u|}{u}du= 2\int_0^{\pi} \sin u (\sum_{k=1}^{n-1}
\frac{1}{u+k\pi})du. 
\]
For $0\leq u\leq \pi,$ the sum is contained between $\frac{1}{\pi}
\sum_{k=2}^n\frac{1}{k}$ and $\frac{1}{\pi} \sum_{k=1}^{n-1}\frac{1}{k},$
and so is strictly of order $\frac{1}{\pi} \ln n.$ Collecting estimates, we
obtain $||D_n||_1= \frac{4}{\pi^2}\ln n+O(1).$

\begin{thm}
There is a continuous function whose Fourier series diverges at a point.
\end{thm}

{\bf Proof:}\quad Suppose it were true that $S_n(h)=\sum_{-n}^n a_k(h)$ has
a limit (i.e. $S_n(h)(e^{ix})$ converges at $x=0$) for any $h\in C(T).$ Then
we would have $|S_n(h)|\leq M_h$ for all $n$, where $M_h$ is a constant. For
each $n$, $S_n$ is a linear functional on $C(T),$ given as $S_n(h)=\int
h(e^{ix})D_n(e^{ix})d\sigma(x).$ If $S_n(h)\leq M_h$ for all $n$ and for
every $h\in C(T),$ then by the Banach-Steinhaus theorem the sequence of
norms $||S_n||$ is bounded. But this norm is the Lebesgue constant, which
tends to $\infty$ with $n.$ This proves that $S_n(h)$ is unbounded for some $
h,$ and thus divergent.

\vspace{0.2in}

{\bf The Fej\'{e}r kernel}

Define 
\[
K_n(e^{ix})=\sum_{-n}^n (1-\frac{|k|}{n})e^{ikx}. 
\]
Then 
\begin{eqnarray*}
f*K_n(e^{ix})&=&\frac{1}{n} \int
f(e^{i(x-t)})\sum_{k=-n}^n(n-|k|)e^{ikt}d\sigma(t) \\
&=&\frac{1}{n}\int f(e^{i(x-t)})(D_{n-1}(e^{it})+
D_{n-2}(e^{it})+\cdots+D_1(e^{it})+D_0)d\sigma(t) \\
&=&\frac{1}{n}\sum_{j=0}^{n-1}(f*D_j)(e^{ix}) =\frac{1}{n}
\sum_{j=0}^{n-1}S_j(f)(e^{ix}).
\end{eqnarray*}

One can easily verify, by multiplying out, that 
\[
\frac{1}{n}|\sum_{k=0}^{n-1} e^{ikx}|^2= \sum_{k=-n}^n (1-\frac{|k|}{n}
)e^{ikx}. 
\]
When we sum the geometric series and simplify, we find 
\[
K_n(e^{ix})= \frac{1}{n}(\frac{\sin \frac{1}{2}nx}{\sin \frac{1}{2}x})^2. 
\]
Thus the Dirichlet and Fej\'{e}r kernels are related by the formula 
\[
K_{2n+1}(e^{ix})= \frac{1}{2n+1}D^2_n(e^{ix}). 
\]

Note that $K_n$ is an approximate identity on $T.$ Thus for any $f\in
L^1(T), $ $K_n*f(e^{ix})\rightarrow f(e^{ix})$ at every point of of
continuity of $f,$ and the convergence is uniform over every closed interval
of continuity. In particular, $K_n*f$ tends to $f$ uniformly everywhere if $
f $ is continuous everywhere. It holds also that if $f\in L^p, 1\leq
p<\infty,$ then $||K_n*f-f||_p\rightarrow 0.$

The functions $K_n$ are trigonometric polynomials; this fact has interesting
consequence.

\begin{enumerate}
\item Since $K_n$'s are infinitely differentiable, any continuous function $
h $ is approximated uniformly by the infinitely differentiable functions (in
fact, trigonometric polynomials) $K_n*h.$

\item We also obtain another proof of the Unicity theorem in $L^1(T).$
Suppose that $a_k(f)=0$ for all $k.$ Then for each $n,$ $
a_k(K_n*f)=a_k(K_n)a_n(f)=0, \forall k.$ Thus the trigonometric polynomial $
K_n*f\equiv 0.$ Since $||K_n*f-f||_1\rightarrow 0,$ $f=0$ a.e.
\end{enumerate}

\vspace{0.1in}

{\bf The Poisson kernel}

Define, for $0<r<1,$ 
\[
P_r(e^{it})=\sum_{-\infty}^{\infty}r^{|n|}e^{int}. 
\]
The series converges absolutely, and we can easily obtain that, if $
z=re^{i\theta}, 0\leq r<1, $ then 
\[
P_r(e^{i(\theta-t)})=Re(\frac{e^{it}+z}{e^{it}-z})= \frac{1-r^2}{1-2r\cos
(\theta-t)+r^2}. 
\]
One can verify that $P_r$ is an approximate identity (with continuous
parameter $r$). Clearly, if $f\in L^1(T)$ then 
\[
\sum_{-\infty}^\infty a_n(f)r^{|n|}e^{in\theta}= =(P_r*f)(e^{i\theta})=\int
P_r(e^{i(\theta-t)})f(e^{it})d\sigma(t).
\]

\begin{thm}
The Poisson integral $(P_r*f)(e^{i\theta})$ provides the harmonic extension
of $f\in L^1(T)$ to the interior of the circle. Moreover,

\begin{enumerate}
\item If $f\in L^p(T)$ with $1\leq p<\infty,$ then $||P_r*f-f||_p\rightarrow
0 $ as $r\uparrow 1.$

\item Let $f\in L^1(T).$ At every point $t$ where $f$ is the derivative of
its integral (hence almost everywhere)  $P_r*f(e^{it})\rightarrow f(e^{it})$
as $r\uparrow 1$ (radial limit). Actually, for almost all $t$, $
P_r*f(e^{i\theta})\rightarrow f(e^{it})$ as $re^{i\theta}\rightarrow e^{it}$
nontangentially. This result depends on particular properties of the Poisson
kernel, and is not true for all approximate identities.
\end{enumerate}
\end{thm}

{\bf Proof:}\quad We prove that $P_r*f(e^{ix})$ is harmonic in $D$ (open
unit disk). If $f$ is real, then $P_r*f$ is the
real part of 
\[
\int \frac{e^{it}+z}{e^{it}-z}f(e^{it})d\sigma(t), 
\]
which is an analytic function of $z=re^{i\theta}$ in $D.$ Hence $P_r*f(e^{i\theta})$ is harmonic in $D.$ Since linear combinations of
harmonic functions are harmonic, $P_r*f(e^{i\theta})$ is a complex harmonic function on $D$ for any $f\in L^1(T),$
the class of all complex, Lebesgue integrable functions on $T.$

\begin{thm}
Suppose $f\in L^1(T)$ and $f\geq 0.$ Then $f$ is the boundary function of a
nonnegative harmonic function. If $f$ is bounded, it is the boundary
function of a harmonic function with the same bounds.
\end{thm}

{\bf Proof:}\quad Let $F(z)=P_r*f(e^{ix}).$ Then $F(z)$ is harmonic in $D$
such that $\lim_{r\uparrow 1}F(re^{i\theta})=f(e^{i\theta})$ for a.e. $
\theta.$ Since $P_r$ is nonnegative, $F(z)$ is certainly positive whenever $f $ is nonnegative. If $|f|\leq M,$ then $||P_r*f||_\infty\leq M ||P_r||_1
=M. $

\begin{thm}
A harmonic function $F$ in $D$ (open disk) is bounded if and only if it is
the Possion integral of some bounded function $f$ on $T.$
\end{thm}

{\bf Proof:}\quad We need only to show the necessity. Let $F$ be harmonic
and bounded in $D.$ Let $r_n\uparrow 1$ and write $f_n(e^{it})=F(r_ne^{it}).$
The sequence $f_n$ is a bounded sequence in $L^\infty(T);$ hence for
some sequence $n_j\rightarrow\infty,$ $f_{n_j}$ converges in the weak-*
topology ($L^\infty(T)$ being the dual of $L^1(T)$) to some function $
f(e^{it}).$

Let $r e^{i\theta}\in D,$ then 
\begin{eqnarray*}
\int f(e^{it})P_r (e^{i(\theta-t)})d\sigma(t) &=&\lim_{j\rightarrow\infty}
\int F(r_{n_j}e^{it})P_r (e^{i(\theta-t)})d\sigma(t) \\
&=&\lim_{j\rightarrow\infty} F(r_{n_j}re^{it})=F(re^{i\theta}).
\end{eqnarray*}

Note that in the above derivation we used the fact that for any harmonic
function $u$ in $D$ and $0\leq r<1,$ 
\[
\int P_r(e^{i(\theta-t)})u(\rho e^{it})d\sigma(t) =u(\rho re^{i\theta}). 
\]
This can be verified by considering the representation theorem of harmonic
functions in disk : If $u$ is real, continuous on $|z|\leq\rho$ and harmonic
in $|z|<\rho,$ then for $z=\rho_1e^{i\theta},$ $\rho_1<\rho,$ 
\[
u(\rho_1e^{i\theta})=\frac{1}{2\pi}\int_0^{2\pi} Re[\frac{\rho e^{it}+z}{
\rho e^{it}-z}]u(\rho e^{it})dt. 
\]
Let $\rho_1=r\rho$ (Note $0\leq r<1$). Then 
\[
u(r\rho e^{i\theta})=\frac{1}{2\pi}\int_0^{2\pi} Re[\frac{e^{it}+r
e^{i\theta}} {e^{it}-re^{i\theta}}]u(\rho e^{it})dt =\frac{1}{2\pi}
\int_0^{2\pi}P_r(e^{i(\theta-t)})u(\rho e^{it})dt. 
\]
For complex harmonic, we consider it as a sum of real part and imaginary
part.

\section{Summability; Metric Theorems}

We have shown the following theorems in the last section:

\begin{thm}
The Poisson integral $(P_r*f)(e^{i\theta})$ provides the harmonic extension
of $f\in L^1(T)$ to the interior of the circle so that

\begin{enumerate}
\item If $f\in C(T),$ then $P_r*f$ converges to $f$ uniformly as $
r\rightarrow 1.$

\item If $f\in L^p(T)$ with $1\leq p<\infty,$ then $||P_r*f-f||_p\rightarrow
0 $ as $r\uparrow 1.$

\item Let $f\in L^1(T).$ At every point $t$ where $f$ is the derivative of
its integral (hence almost everywhere) $P_r*f(e^{i\theta})\rightarrow
f(e^{it})$ as $re^{i\theta}\rightarrow e^{it}$ nontangentially.
\end{enumerate}
\end{thm}

\begin{thm}
Let $\mu$ be any finite complex measure on $T$. Then 
\[
P_r*\mu(e^{i\theta})=\int P_r(e^{i(\theta-t)})d\mu(t) 
\]
is a harmonic function in $D$ and converges to $\mu$ in the weak* topology
of $M(T).$ That is, for any $h\in C(T),$ 
\[
\int h(e^{ix})(P_r*\mu)(e^{ix})d\sigma(x)\rightarrow \int h(e^{ix})d\mu(x),
\quad r\rightarrow 1. 
\]
\end{thm}

{\bf Proof:}\quad $P_r*\mu(e^{i\theta})$ is a continuous function of $\theta$
and $a_n(P_r*\mu)=a_n(P_r)a_n(\mu).$ Thus 
\[
P_r*\mu(e^{i\theta})=\sum a_n(\mu)r^{|n|}e^{in\theta}. 
\]
Since $|a_n(\mu)|\leq |\mu|(T),$ the above series converges uniformly in
every compact subset of $D$ and so $P_r*\mu(e^{i\theta})$ is harmonic in $D.$

For any $h\in C(T),$ 
\[
\int h(e^{ix})\int (P_r(e^{i(x-t)})d\mu(t)d\sigma(x)= \int \int
h(e^{ix})P_r(e^{i(x-t)})d\sigma(x)d\mu(t). 
\]
Note that $P_r*h$ converges uniformly to $h.$ Taking the limit of above
equality as $r\uparrow 1,$ we complete the proof of the theorem.

\vspace{0.5in}

\begin{thm}
A harmonic function $F$ in $D$ (open disk) is bounded if and only if it is
the Possion integral of some bounded function $f$ on $T.$
\end{thm}

{\bf Proof:}\quad We need only to show the necessity. Let $F$ be harmonic
and bounded in $D.$ Let $r_n\uparrow 1$ and write $f_n(e^{it})=F(r_ne^{it}).$
The sequence $f_n$ is a bounded sequence in $L^\infty(T);$ hence for
some sequence $n_j\rightarrow\infty,$ $f_{n_j}$ converges in the weak-*
topology ($L^\infty(T)$ is the dual of $L^1(T)$) to some function $f(e^{it}).$

Let $r e^{i\theta}\in D,$ then 
\begin{eqnarray*}
\int f(e^{it})P_r (e^{i(\theta-t)})d\sigma(t) &=&\lim_{j\rightarrow\infty}
\int F(r_{n_j}e^{it})P_r (e^{i(\theta-t)})d\sigma(t) \\
&=&\lim_{j\rightarrow\infty} F(r_{n_j}re^{it})=F(re^{i\theta}).
\end{eqnarray*}
Clearly, the fact that $F=P_r*f$ implies immediately that if $F(z)=\sum
a_nr^{|n|}e^{in\theta},$ then $a_n(f)=a_n.$

Note that in the above derivation we used the fact that for any harmonic
function $u$ in $D$ and $0\leq r<1,$ 
\[
\int P_r(e^{i(\theta-t)})u(\rho e^{it})d\sigma(t) =u(\rho re^{i\theta}). 
\]
This can be verified by considering the representation theorem of harmonic
functions in disk: If $u$ is real, continuous on $|z|\leq\rho$ and harmonic
in $|z|<\rho,$ then for $z=\rho_1e^{i\theta},$ $\rho_1<\rho,$ 
\[
u(\rho_1e^{i\theta})=\frac{1}{2\pi}\int_0^{2\pi} Re[\frac{\rho e^{it}+z}{
\rho e^{it}-z}]u(\rho e^{it})dt. 
\]
Let $\rho_1=r\rho$ ($0\leq r<1$). Then 
\[
u(r\rho e^{i\theta})=\frac{1}{2\pi}\int_0^{2\pi} Re[\frac{e^{it}+r
e^{i\theta}} {e^{it}-re^{i\theta}}]u(\rho e^{it})dt =\frac{1}{2\pi}
\int_0^{2\pi}P_r(e^{i(\theta-t)})u(\rho e^{it})dt. 
\]
For complex harmonic, we consider it as a sum of real part and imaginary
part.

\vspace{0.1in}

Let $C(T)$ be the space of continuous functions on $T.$ Then by F. Riesz'
theorem, $(C(T))^*=M(T),$ where $M(T)$ is the space of bounded complex Borel
measures on $T.$ Since $C(T)$ is separable, every bounded subset of $M(T)$
is weak* sequentially compact. Note that $L^1(T)$ is contained in $M(T),$ if
we identify $f\in L^1$ with the measure $f(e^{ix})d\sigma(x).$

\begin{lem}
Every complex harmonic function $F$ in $D$ has a development $
F(re^{i\theta})=\sum a_nr^{|n|}e^{in\theta}.$
\end{lem}

{\bf Proof:} Suppose that $F$ is real. Since $D$ is a simply connected
region, $F$ has a harmonic conjugate $G$ so that $H=F+iG$ is analytic in $D.$
We write $H(z)=\sum_{n=0}^\infty a_nz^n.$ Then 
\begin{eqnarray*}
H(z)=Re(H)&=&a_0+\frac{1}{2}\biggr(\sum_{n=1}^\infty a_nr^ne^{in\theta}
+\sum_{n=1}^\infty \overline{a_n}r^ne^{-inx}\biggl) \\
&=&a_0+\frac{1}{2}\sum_{n=-\infty}^\infty a_nr^{|n|}e^{in\theta}
\end{eqnarray*}
where $a_{-n}=\overline{a_{n}}$ for $n=1,2,\cdots.$ If $F$ is complex, then
it is linear combination of two real harmonic functions. The desired
development of $F$ is the sum of two absolutely convergent series.

\begin{thm}
Let $F$ be a function harmonic in the unit disk $D.$ There is a unique
measure $\mu\in M(T)$ such that 
\[
F(z)=P_r*\mu(e^{i\theta}),\qquad z=re^{i\theta} 
\]
if and only if 
\[
A_r=\int |F(re^{i\theta})|d\sigma(\theta)\leq K,\quad \forall\quad 0<r<1. 
\]
Moreover, $||\mu||=\lim_{r\uparrow 1}A_r.$
\end{thm}

{\bf Proof:} Necessity: If we think of $P_r*\mu(e^{i
\theta}) $ as a family (with continuous parameter $0<r<1$) of functions
defined on $T$, call them $f_r,$ then 
\[
||f_r||_1=||P_r*\mu||\leq ||P_r||_1||\mu||. 
\]
Therefore, $\{f_r\}$ in bounded in $L^1(T).$ This shows that in order for a
harmonic function $F(z)$ in $D$ to be the Poisson integral of some measure
it is necessary that 
\[
\int |F(re^{i\theta})|d\sigma(\theta)\leq K,\quad \forall\quad 0<r<1. 
\]

Sufficiency: Given a harmonic function $F(re^{ix})$ in $D,$ by
the lemma, one can always write $F$ as a `power' series: 
\[
F(re^{ix})=\sum_{n=-\infty}^\infty a_nr^{|n|}e^{inx}, 
\]
where $a_nr^{|n|}$ is the $n$th Fourier coefficient of $
f_r(e^{ix})=F(re^{ix}).$

By assumption, $||f_r||_1\leq K,$ i.e., $||f_rd\sigma||_{M(T)}=||f_r||_1\leq
K \quad \forall r.$ Since $C(T),$ as the predual of $M(T),$ is separable
normed space (polynomials with rational coefficients are dense in $C(T)$),
by Banach-Alaoglu theorem, the closure of $\{f_rd\sigma\}$ in $M(T)$ is weak*
sequentially compact. Therefore, there is a subsequence $\{f_{r_j}\}d\sigma$
of $f_rd\sigma$ that converges to some $\mu\in M(T)$ in weak* topology. That
is, 
\[
\int h(e^{-ix})f_{r_j}(e^{ix})d\sigma(x)\rightarrow \int
h(e^{-ix})d\mu,\qquad j\rightarrow\infty 
\]
for each $h\in C(T).$ In particular, for each $n,$ $a_n(f_{r_j})\rightarrow
a_n(\mu)$ as $j\rightarrow\infty.$ On the other hand, $a_n(f_r)=a_nr^{|n|}
\rightarrow a_n$ as $r\uparrow 1.$ Therefore, $a_n(\mu)=a_n$ for all $n.$

It follows from the Unicity theorem that $\mu$ is uniquely determined by $
a_n,$ therefore by $F,$ and that since $a_n(f_r)=a_nr^{|n|}=a_n{\mu}
r^{|n|}=a_n(P_r*\mu),$ $f_r=P_r*\mu,$ i.e. $F(re^{ix})=P_r*\mu(e^{ix}).$

We show that $||\mu||=\lim_{r\uparrow 1}A_r.$ Note that $\mu=\lim_{j
\rightarrow \infty}F(r_je^{ix})d\sigma(x)$ in the weak* topology of $M(T)$
as the dual of $C(T).$ It follows that $||\mu||\leq \liminf_{j\rightarrow
\infty} A_{r_j}$ where $A_{r_j}=||F(r_je^{i\cdot})||_1.$ Since $A_r$ increases with $r$ and $A_r\leq K,$ $||\mu||\leq
\lim_{r\rightarrow\infty}A_r.$ Furthermore, the inequality cannot be strict.
Note that $f_r=P_r*\mu$ and $||f_r||_1\leq ||P_r||_1||\mu||.$ Therefore, $
A_r=||f_r||_1\leq ||\mu||$ for every $0<r<1.$ If the inequality were strict,
we would have $A_r\leq ||\mu||<\lim_{r\rightarrow\infty}A_r$ for $0<r<1,$
which is impossible.

As to the norm convergence of $||f_r-\mu||_{M(T)}\rightarrow 0$ as $
r\rightarrow 0,$ if $\mu$ is absolutely continuous, then $\mu=fd\sigma$ for
some $f\in L^1(T).$ Now $f_r=P_r*\mu$ is indeed $f_r=P_r*f.$ Thus, by
Fej\'{e}r's theorem, $||f_r-f||_1\rightarrow 0.$ That is, $||f_r-\mu||_{M(T)}
\rightarrow 0.$ \hfill 

\vspace{0.1in}

\begin{thm}
Suppose $F$ is harmonic in the unit disk $D.$ If $1<p\leq\infty,$ $F$ is the
Possion integral of a (unique) function $f$ belonging to $L^p(T)$ if and
only if 
\[
A_r^p=\int |F(re^{ix})|^pd\sigma(x)\leq K<\infty,\quad 0<r<1. 
\]
Moreover, $||f||_p =\lim_{r\uparrow 1}A_r.$
\end{thm}

{\bf Proof:}\quad Similar to the proof of the previous theorem.

\vspace{0.1in}

\section{Convergence a.e.}

\begin{thm}
If $f\in L^1(T),$ then the symmetric derivative of indefinite integral of $f$ 
\[
\lim_{t\rightarrow 0}\frac{1}{2t}\int_{x-t}^{x+t}f(u)du=f(x),\quad a.e.
\]
Furthermore, for each $x$ where this holds, 
\[
\lim_{r\uparrow 1}(p_r*f)(x)=f(x). 
\]
\end{thm}

{\bf Proof:}\quad The first part is simply the Lebesgue differentiation
Theorem: if $f\in L^1,$ then for a.e. x, 
\[
\lim_{Q\downarrow x}\frac{1}{|Q|} \int_Qf(u)du=f(x), 
\]
where $Q$'s are intervals centered at $x.$

Let $x$ be a point where the above limit holds, i.e., at $x$ 
\[
\int_{0}^t(f(x+u)+f(x-u)-2f(x))du=o(t). 
\]
Let $G(t)=\int_{0}^t(f(x+u)+f(x-u)-2f(x))du.$ Given any $\epsilon>0,$ there
is $\delta>0$ so that if $0\leq t\leq\delta,$ then $|G(t)|\leq\epsilon t.$
Also note that $G(t)$ is absolutely continuous so that $G^{\prime}(t)=
f(x+t)+f(x-t)-2f(x)$ for a.e. $t$. We write 
\begin{eqnarray*}
(p_r*f)(x)&=&\frac{1}{2\pi}\int_0^\pi (f(x+u)+f(x-u)-2f(x))p_r(u)du \\
&=&\frac{1}{2\pi}\biggr(\int_0^\delta+ \int_\delta^\pi\biggl)
(f(x+u)+f(x-u)-2f(x))p_r(u)du=I_1+I_2.
\end{eqnarray*}

We consider $I_1$ first. We have 
\begin{eqnarray*}
2\pi|I_1|&=&|\int_0^\delta G^{\prime}(u)p_r(u)du| \\
&=&|G(\delta)p_r(\delta)-\int_0^\delta G(u)p_r^{\prime}(u)du| \\
&\leq&|G(\delta)p_r(\delta)|+\int_0^\delta |G(u)||p_r^{\prime}(u)|du \\
&\leq&\epsilon\biggl(|\delta p_r(\delta)|+\int_0^\delta u|p_r^{\prime}(u)|du
\biggr). \qquad
\end{eqnarray*}

Note that $\delta p_r(\delta)>0.$ Using integration by parts, we have 
\[
|\delta p_r(\delta)|=\delta p_r(\delta)= \int_0^\delta
up_r^{\prime}(u)du+\int_0^\delta p_r(u)du. 
\]
Furthermore, since $p_r^{\prime}(u)$ keeps constant sign (negative) as $x>0,$
\[2\pi |I_1|= \epsilon (\int_0^\delta up_r^{\prime}(u)du+\int_0^\delta
p_r(u)du +\int_0^\delta u|p_r^{\prime}(u)|du)=\epsilon\int_0^\delta
p_r(u)du. \]
The last integral goes to $0$ as $r\uparrow 1.$

\vspace{0.1in}

The following lemma can be used to give an alternative proof of the above
theorem:

\begin{lem}
Let $e_n(u)$ be an even approximate identity. Then for any bounded function $f,$ at the point $x$ where $f(x+)$ and $f(x-)$ exist, we have 
\[\lim_{n\rightarrow \infty}\int_{-\pi}^{\pi}
f(x-u)e_n(u)d\sigma(u)\rightarrow \frac{f(x+)+f(x-)}{2}. \]
In particular, if $\lim_{u\rightarrow x} f(u)=L,$ then 
\[\lim_{n\rightarrow \infty}\int_{-\pi}^{\pi}
f(x-u)e_n(u)d\sigma(u)\rightarrow L. \]
\end{lem}

{\bf Proof:}\quad If $e_n$ is even, then 
\begin{eqnarray*}
I&=&\int_{-\pi}^{\pi}f(x-u)e_n(u)d\sigma(u) -\frac{f(x+)+f(x-)}{2} \\
&=&\int_0^{\pi}\biggl(f(x+u)+f(x-u)-f(x+)-f(x-)\biggr)e_n(u)d\sigma(u).
\end{eqnarray*}
Given $\epsilon>0$, there exists $\delta>0$ such that for $0< u\leq\delta,$ $
|f(x+u)-f(x+)|<\epsilon$ and $|f(x-u)-f(x-)|<\epsilon.$ We write 
\[
I=\biggl(\int_0^{\delta}+ \int_\delta^\pi\biggr)
(f(x+u)-f(x+)+f(x-u)-f(x-))e_n(u)d\sigma(u)=I_1+I_2. 
\]
For $I_1,$ we have 
\[
|I_1|\leq 2\epsilon\int_0^\pi e_n(u)d\sigma(u)=2\epsilon. 
\]
For $I_2,$ we have 
\[
|I_2|\leq 4M\int_\delta^\pi e_n(u)d\sigma(u)<\epsilon, 
\]
for sufficiently large $n.$ 

\vspace{0.2in}

An alternative proof of the above theorem.

{\bf Proof:}\quad First, we make the following assumptions successively:

\begin{enumerate}
\item We may assume that $x=0$ is the point where 
\[
\lim_{t\rightarrow 0}\frac{1}{2t}\int_{x-t}^{x+t}f(u)du=f(x). 
\]
That is, 
\[
\lim_{t\rightarrow 0}\frac{1}{2t}\int_{-t}^{t}f(u)du=f(0). 
\]
Assume that the limit holds for $f$ at $x=a.$ Then $g(x)=f(x+a)$ satisfies 
\[
\lim_{t\rightarrow 0}\frac{1}{2t}\int_{-t}^{t}g(u)du=g(0). 
\]
If the theorem is proved for $g$ at $0,$ then $(p_r*g)(0)\rightarrow g(0)$ is
simply $(p_r*f)(a)\rightarrow f(a).$

\item We may also assume that $f(0)=0.$ Let $g(x)=f(x)-f(0).$ Then $g(0)=0.$
If the theorem is proved for $g,$ then $(p_r*g)(0)\rightarrow 0)$ is simply $
(P_r*(f(\cdot)-f(0))(0) =(p_r*f)(0)-f(0)\rightarrow 0,$ which is $
(p_r*f)(0)\rightarrow f(0).$

\item Finally, we may assume that ${\displaystyle\int f(x)d\sigma(x)=0}.$
Let $g$ be a smooth function with ${\displaystyle\int g=\int f}$, and
vanishing on a neighborhood of $x=0$ (maintaining the above two
assumptions). If the theorem is proved for $f-g$, then it follows for $f.$
\end{enumerate}

Now we proceed to prove the theorem:  If $f\in L^1(T)$ with 
${\displaystyle\int f}=0,$ and 
\[
\lim_{t\rightarrow 0}\frac{1}{2t}\int_{-t}^{t}f(u)du=0, 
\]
then 
\[
\lim_{r\uparrow 1}(p_r*f)(0)=0. 
\]

{\bf Proof}: \quad We first prove that 
\[
q_r(x)=r^{-1}(-\sin x p_r^{\prime}(x)) 
\]
is an approximate identity. Furthermore, if we define $F(t)={\displaystyle
\int_{-\pi}^tf(x)dx},$ then we can write $p_r*f(0)$ as, using integration by
parts and ${\displaystyle\int f}=0,$ 
\[
p_r*f(0)=\frac{1}{2\pi}\int_{-\pi}^\pi p_r(x)f(x)dx=-\frac{1}{2\pi}
\int_{-\pi}^\pi p_r^{\prime}(x)F(x)dx. 
\]
Since $p^{\prime}_r$ is odd, the last integral can be written as 
\begin{eqnarray*}
-\frac{1}{2\pi}\int_{-\pi}^\pi p_r^{\prime}(x)F(x)dx &=& -\frac{1}{2\pi}
\int_{0}^\pi p_r^{\prime}(x)(F(x)-F(-x))dx \\
&=& \frac{1}{\pi}\int_{0}^\pi rq_r(x)\frac{x}{\sin x} \frac{F(x)-F(-x)}{2x}
dx.
\end{eqnarray*}
Note that the function ${\displaystyle\frac{x}{\sin x}\frac{F(x)-F(-x)}{2x}}$
is bounded (since $f\in L^1$, $F$ is bounded. In addition, ${\displaystyle
\frac{F(x)-F(-x)}{2x}\rightarrow 0}$ as $x\rightarrow 0$). By the lemma, we
have that 
\begin{eqnarray*}
p_r*f(0)&=&\frac{1}{\pi}\int_{0}^\pi rq_r(x)\frac{x}{\sin x} \frac{F(x)-F(-x)
}{2x}dx \\
&\rightarrow& \lim_{x\rightarrow 0} (\frac{x}{\sin x}\frac{F(x)-F(-x)}{2x}
)=0, \quad r\uparrow 1.
\end{eqnarray*}

\begin{lem}
The Fej\'{e}r kernel 
\[
K_n(x)=\frac{1}{n}|\sum_{0}^{n-1}e^{ikx}|^2 =\frac{1}{n}(\frac{sin\frac{1}{2}
nx}{sin\frac{1}{2}x})^2 
\]
has a bell-shaped majorant 
\[
K^*_n(x)=\frac{2\pi^2 n}{1+n^2x^2}. 
\]
\end{lem}

{\bf Proof:}\quad The first formula for $K_n$ gives $K_n(x)\leq n$ (used for
smaller $x$). By Jordan's inequality, the second formula leads to $K_n(x)\leq
\frac{\pi^2}{nx^2}$ (used for large $x$). Combining these two gives $
K_n(x)\leq K_n^*(x),$ where $K^*_n(x)=\frac{2\pi^2n}{1+n^2x^2}$ (consider $
|x|\leq \frac{1}{n}$ and $|x|>\frac{1}{n}$ separately).

\begin{thm}
Let $f\in L^1$ and let $x$ be such a point where 
\[
\lim_{t\rightarrow 0}\frac{1}{t}\int_0^t|f(x+u)+f(x-u)-2L|du=0 
\]
for some $L$ (Note that if $f\in L^1,$ then $L=f(x)$ for a.e. x). Then $
K_n*f(x)\rightarrow L$ as $n\rightarrow\infty.$
\end{thm}

{\bf Proof:}\quad First, we write 
\begin{eqnarray*}
|K_n*f(x)-L|&=&|\frac{1}{2\pi}\int_0^\pi K_n(t)(f(x+t)+f(x-t)-2L)dt| \\
&=&|\frac{1}{2\pi}\biggl(\int_0^{\frac{2\pi}{n}}+ \int_{\frac{2\pi}{n}}^\pi
\biggr) K_n(t)(f(x+t)+f(x-t)-2L)dt| \\
&\leq&\frac{1}{2\pi}\biggl(\int_0^{\frac{2\pi}{n}}+ \int_{\frac{2\pi}{n}}^\pi
\biggr) K_n(t)|f(x+t)+f(x-t)-2L|dt=I_1+I_2.
\end{eqnarray*}
It is interesting to note that the splitting point varies with $n.$

Let $H(t)=\int_0^t|f(x+t)+f(x-t)-2L|dt.$ Then $\frac{H(t)}{t}\rightarrow 0$
as $t\rightarrow 0.$ For $I_1$ we have 
\[
2\pi I_1\leq nH(\frac{2\pi}{n})\rightarrow 0,\qquad n\rightarrow \infty. 
\]
For $I_2$ we have 
\begin{eqnarray*}
2\pi I_2&\leq& \int_{\frac{2\pi}{n}}^\pi K_n^*(t) |f(x+t)+f(x-t)-2L|dt \\
&=&K^*_n(\pi)H(\pi)-K^*_n(\frac{2\pi}{n})H(\frac{2\pi}{n}) -\int_{\frac{2\pi
}{n}}^\pi H(t)(K_n^*)^{\prime}(t) dt.
\end{eqnarray*}
Since $H(\pi)$ is finite, $K^*_n(\pi)H(\pi)\rightarrow 0$ as $n\rightarrow
0. $ We also have that 
\[
K^*_n(\frac{2\pi}{n})H(\frac{2\pi}{n})\leq \frac{n}{4}H(\frac{2\pi}{n})
\rightarrow 0,\quad n\rightarrow\infty. 
\]
As for the last term, we have 
\begin{eqnarray*}
&-&\int_{\frac{2\pi}{n}}^\pi H(t)(K_n^*)^{\prime}(t) dt \\
&=&\int_{\frac{2\pi}{n}}^\pi H(t)\frac{4\pi^2 n^3t}{(1+n^2t^2)^2}dt \\
&\leq&\int_{\frac{2\pi}{n}}^\pi H(t)\frac{4\pi^2}{nt^3}dt.
\end{eqnarray*}
We show that the last integral tends to $0.$ Given $\epsilon>0$, there
exists $\delta>0$ such that $H(t)<\epsilon t$ if $0<t<\delta.$ Then 
\[
\frac{1}{n}\int_{\frac{2\pi}{n}}^\pi \frac{H(t)}{t^3}dt \leq \frac{1}{n}
\int_{\frac{2\pi}{n}}^\delta \frac{\epsilon }{t^2}dt +\frac{1}{n}
\int_{\delta}^\pi \frac{H(t)}{t^3}dt. 
\]
The first term on the right is majorized by $\epsilon$ (simply integrate),
while the last term clearly tends to $0$ as $n\rightarrow \infty.$ 

\begin{thm}
The boundedness of $f$ in the lemma 6.1 is indispensable.
\end{thm}

{\bf Proof:}\quad We define $e_n$ as an (even) approximate identity with a
``tail'': 
\[
2\pi e_n(x)=\left\{
\begin{array}{ll}
n-\sqrt{n} & 0\leq x<\frac{1}{n}, \\ 
0 & \frac{1}{n}\leq x<(\pi-\frac{1}{n}), \\ 
\sqrt{n} & (\pi-\frac{1}{n})\leq x<\pi.
\end{array}
\right. 
\]
Clearly, $e_n(x)\geq 0$, ${\displaystyle \int e_n(x)dx=1},$ and for any $
\delta>0,$ $\int_{|x|>\delta}e_n(x)dx =\frac{1}{\sqrt{n}}\rightarrow 0.$

Let $f(x)$ be an even function defined as 
\[
f(x)=\left\{
\begin{array}{ll}
0 & 0\leq x<\pi-1 \\ 
\sqrt{n} & \pi-\frac{1}{n}\leq x<\pi-\frac{1}{n+1}, n=1,2,\cdots.
\end{array}
\right. 
\]
Then ${\displaystyle \int f=\sum \frac{\sqrt{n}}{n(n+1)}<\infty}$ so that $
f\in L^1.$ But $e_n*f(0)=\int e_n(x)f(x)dx\geq 1$ for all $n.$ We do not
have $e_n*f(0)\rightarrow f(0)=0.$ In this example, the coincidence of the
tail of $f$ and that of $e_n$ makes ${\displaystyle \int e_n(x)f(x)dx}$ big.

\begin{thm}
A bounded analytic function $F$ in $D$ has radial limits at a.e. boundary
point.
\end{thm}

{\bf Proof:}\quad We write $F=G+iH.$ Since $G$ is harmonic and bounded,
there is a bounded function $g$ on $T$ such that $
G(re^{ix})=(p_r*g)(e^{ix}). $ Note that $(p_r*g)(e^{ix}) \rightarrow
g(e^{ix}), a.e.$ That is, $G(re^{ix})\rightarrow g(e^{ix}), a.e.,$ as $
r\uparrow 1.$ Similarly, there is a bounded function $h$ on $T$ so that $
H(re^{ix})\rightarrow h(e^{ix}).$ Thus, $F(re^{ix})\rightarrow
g(e^{ix})+ih(e^{ix})$ as $r\uparrow 1.$ 

\begin{thm}
If $\mu$ is a singular measure on $T,$ then $P_r*\mu$ tends to $0$ a.e. as $
r $ increases to $1$.
\end{thm}

{\bf Proof:}\quad Let $x=0$ be a point where 
\[
\lim_{u\rightarrow 0} \frac{\mu([-u,u))}{2u}=0. 
\]
We show that 
\[
(p_r*\mu)(0)=\int_{-\pi}^\pi p_r(u)d\mu(u)\rightarrow 0,\quad r\uparrow 1. 
\]
Let $F(t)={\displaystyle\int_{-\pi}^td\mu(v)}.$ Then, by Fubini's theorem, 
\begin{eqnarray*}
\int_{-\pi}^\pi p^{\prime}_r(u)F(u)du &=&\int_{-\pi}^\pi
p_r^{\prime}(u)\int_{-\pi}^ud\mu(v)du \\
&=&\int_{-\pi}^\pi \int_v^\pi p_r^{\prime}(u)dud\mu(v) \\
&=&\int_{-\pi}^\pi (p_r(\pi)-p_r(v))d\mu(v) \\
&=&p_r(\pi)\mu([-\pi,\pi))+\int_{-\pi}^\pi p_r(v)d\mu(v).
\end{eqnarray*}
It follows that 
\[
\int_{-\pi}^\pi p_r(v)d\mu(v)= \int_{-\pi}^\pi
p^{\prime}_r(u)F(u)du-p_r(\pi)\mu([-\pi,\pi)). 
\]
For the last integral, noting that $p^{\prime}_r$ is odd, we have
\begin{eqnarray*}
\int_{-\pi}^\pi p^{\prime}_r(u)F(u)du &=&\int_0^\pi
p^{\prime}_r(u)(F(u)-F(-u))du \\
&=&\int_0^\pi 2up^{\prime}_r(u)\frac{\mu([-u,u))}{2u}du \\
&\rightarrow& \lim_{u\rightarrow 0}\frac{\mu([-u,u))}{2u}=0,\quad r\uparrow 1.
\end{eqnarray*}
\hfill

\section{Herglotz' Theorem}

\begin{df}
A complex sequence $\{u_n\}_{n=-\infty}^\infty$ is called positive definite
if 
\[
\sum_{m,n}u_{m-n}c_m\overline{c_n}\geq 0 
\]
for every sequence $\{c_n\}_{n=-\infty}^\infty$ such that $c_n=0$ except for
a finite number of $n.$
\end{df}

\begin{thm}
Let $\mu$ be any positive measure on $[0,2\pi).$ Set 
\[
u_n=\int e^{-inx}d\mu(x). 
\]
Then $u_n$ is positive definite.
\end{thm}

{\bf Proof:}\quad Note that 
\begin{eqnarray*}
\int |\sum c_me^{-imx}|^2d\mu(x)&=& \int (\sum c_me^{-imx})(\sum \overline{
c_n} e^{inx})d\mu(x) \\
&=&\int \sum_{m,n} c_m\overline{c_n} e^{-i(m-n)x}d\mu(x) \\
&=& \sum_{m,n}c_m\overline{c_n}u_{m-n}.
\end{eqnarray*}

\vspace{0.1in}

The result to be proved is that these are the only positive definite
sequences.

\begin{thm}[Herglotz]
Every positive definite sequence $u_n$ is the sequence of the
Fourier-Stieltjes coefficients of a positive measure.
\end{thm}

{\bf Proof:}\quad Use the given sequence $u_n$ to define a linear
functional $F$ on trigonometric polynomials. For each trigonometric
polynomial $\varphi =\sum c_ne^{inx},$ define 
\[
F(\varphi)=\sum c_nu_n. 
\]
Clearly, $F$ is a linear (over complex scalar) functional defined on a dense
subspace of $C(T).$

\vspace{0.1in}

Prove that $F$ is nonnegative on every nonnegative trigonometric
polynomial. Note that 
\[
|\varphi|^2=(\sum c_me^{imx})(\sum \overline{c_n}e^{-inx}) =\sum_{m,n}c_m
\overline{c_n}e^{i(m-n)x} 
\]
and that $u_n$ is positive definite. Thus 
\[
F(|\varphi|^2)=F(\varphi \overline{\varphi}) =\sum_{m,n}c_m\overline{c_n}
u_{m-n}\geq 0. 
\]
This means that $F(\psi)\geq 0$ for any trigonometric polynomial $\psi$ that
can be written as $\psi=|\varphi|^2$ for some trigonometric polynomial $
\varphi.$ Note that every nonnegative trigonometric polynomial has this
form. Thus $F$ is a linear functional on
trigonometric polynomials that is nonnegative on every nonnegative
trigonometric polynomial.

\vspace{0.1in}

Prove that $F$ is a (complex-valued) continuous linear functional on
the space of real trigonometric polynomials. It follows from the above
argument that for any two trigonometric polynomials $\varphi$ and $\psi$
with $\varphi\leq \psi,$ $F(\varphi)\leq F(\psi).$ In particular, if we
denote $k=F(1),$ then (choose $\psi=0$ and $1,$ respectively) for any $0\leq
\varphi\leq 1,$ $0\leq F(\varphi)\leq k.$ Let $\varphi$ be a trigonometric
polynomial with $|\varphi|\leq 1$. Then $0\leq 1+\varphi\leq 2.$ Thus $0\leq
F(1)+F(\varphi)=F(1+\varphi)\leq 2k$ and so $|F(\varphi)|\leq k$. Thus $F$
is a (complex-valued) continuous linear functional on the space of real
trigonometric polynomials.

\vspace{0.1in}

Extend $F$ to a continuous linear functional on $C(T)$ - the space of the
complex valued continuous functions.  Since the space of real trigonometric
polynomials is dense in the space of all real continuous function and $F$ is
continuous on the space of real trigonometric polynomials, $F$ can be
extended (by limit) to a continuous linear functional on all real continuous
functions. Further, for complex valued continuous function $f,$ define $
F(f)=F(Ref)+iF(Imf).$ Then $F$ is a continuous linear functional on $C(T)$.
Also note that this definition, when restricted on trigonometric
polynomials, coincides with the original definition for $F.$

\vspace{0.1in}

Use $F$ to find the desired measure $\mu$. Since $F\in C(T)^*,$ by the
Riesz representation theorem there exists a (complex, in general) measure $\mu$ such that 
\[F(f)=\int fd\mu(x)\]
for all $f\in C(T).$ Moreover, looking at how $\mu$ is constructed, we see
that
\[\mu(E)=\int_E \overline{g}d\iota \]
for some $g\in L_{\infty}(T)$ with $||g||_\infty=1$ and $\iota$ a positive
Borel measure on $T.$

\vspace{0.1in}

Prove that $g$ is positive and so $\mu$ is indeed a positive measure.
Since each nonnegative continuous function is the uniform limit of
nonnegative trigonometric polynomials (for example, the Fej\'{e}r means of its
Fourier series), $F$ is nonnegative on all nonnegative continuous
functions. That is, 
\[
F(f)=\int_X f\overline{g}d\iota\geq 0 
\]
for all nonnegative $f\in C(T).$ It follows that $g\geq 0$ a.e. and so $\mu$
is a positive measure.

\vspace{0.1in}

Prove that $\hat{\mu}(n)=u_n.$ It follows from the previous proofs (let $\varphi=e^{inx}$) that $\hat{\mu}(n)=\int e^{-inx}d\mu(x)=F(e^{inx})=u_n.$

\vspace{0.3in}

{\bf An alternative proof of Herglotz' theorem.}

\vspace{0.1in}

Define $A$ as the space of functions $f\in C(T)$ with absolutely convergent
Fourier series. That is, if $f\sim \sum a_n(f)e^{inx},$ then $\sum |a_n(f)|$
converges (so that the Fourier series of $f$ converges uniformly to $f$).

Claim: $A$ is a Banach algebra under multiplication (the product of $f$
and $g$ is defined as $f\overline{g}$) in the norm inherited from $l^1.$

Proof: Let $f\in A.$ $||f||=||f||_A=||\{a_n(f)\}||_{l^1}.$ Then $A$ is
a normed linear space. We prove that $A$ is complete. Let $f_k$ be a Cauchy
sequence in $A,$ that is, $\{a_n(f_k)\}$ is a Cauchy sequence in $l^1.$
Assume that this sequence converges to $\{b_n\}\in l^1.$ Let $f=\sum
b_ne^{inx}.$ Then $f\in A$ and $||f_k-f||
=||\{a_n(f_k)-b_n\}||_{l^1}\rightarrow 0.$

To prove that $A$ is a Banach algebra, we prove first that if both $f$ and $g $ are in $A,$ then so is $f\overline{g}.$ Note that if $f(e^{ix})=\sum a_k(f)e^{ikx}$ and $g(e^{ix})=\sum
a_j(g)e^{ijx},$ then 
\begin{eqnarray*}
f(e^{ix})\overline{g(e^{ix})}&=& (\sum a_k(f)e^{ikx})(\sum \overline{a_j(g)}
e^{-ijx}) \\
&=&\sum_{k,j}a_k(f)\overline{a_j(g)}e^{i(k-j)x} =\sum c_le^{ilx},
\end{eqnarray*}
where 
\[
c_l=\sum_ka_k(f)\overline{a_{k-l}(g)}. 
\]
If $a_k(f), a_j(g)\in l^1,$ then $c_l$ converges absolutely for every $l.$ 
 Moreover, since $\{c_l\}$ is the
convolution of $\{a_k(f)\}$ and $\{a_j(g)\}\in l^1,$ $\{c_l\}\in l^1$ (like
that in $L^1$) and $fg\in A.$ Secondly, we verify that $||f\overline{g}
||\leq ||f||||g||.$ Note that the inequality actually says that $
||a*b||_{l^1} \leq ||a||_{l^1}||b||_{l^1}$ for $a, b\in l^1.$ But it is true
just like in $L^1.$

Claim: Define for all $\varphi=\sum c_ne^{inx}\in A,$ 
\[
F(\varphi)=\sum c_nu_n. 
\]
Then $F(|\varphi|^2)\geq 0$ for each $\varphi\in A.$

We need to justify this definition first. Note that $|u_n|\leq u_0$ for all $n.$ Thus $\sum c_nu_n$ converges absolutely so that $
F(\varphi) $ is well-defined for all $\varphi\in A.$

Next, we show $F(|\varphi|^2)\geq 0$ for all $\varphi\in A.$ Note that $
|\varphi|^2=\varphi\overline{\varphi}\in A.$ By definition of $F$, 
\[
F(|\varphi|^2)=\sum_{k,j}a_k(\varphi)\overline{a_j(\varphi)}u_{k-j}. 
\]
(The coefficient of $u_0$ is $\sum |a_k(\varphi)|^2$). If $
a_k(\varphi) $'s are zeros except for finitely many $k$, then $
\sum_{k,j}a_k(\varphi)\overline{a_j(\varphi)}u_{k-j}\geq 0.$ Hence $
F(|\varphi|^2)\geq 0$ as long as the sum that evaluates $F(|\varphi|^2)$
converges, which is indeed the case because $u_n$'s are bounded.

Claim: If $\psi\in A$ and $\psi>0$ (strictly positive!), then $
\psi=|\varphi|^2$ for some $\varphi\in A.$ Hence, $F(\psi)\geq 0$ for all $\psi\in A$ and $\psi>0.$

Proof: Note that $\sqrt{z}$ is analytic on the right-half (open) plane
that contains the range of $\psi$ ($\psi>0$) and $\psi$ has absolutely
convergent Fourier series. By the Weiner-Levy theorem, $\sqrt{\psi}$ has
absolutely convergent Fourier series. That is, $\sqrt{\psi}\in A.$

Claim: $|F(\varphi)|\leq M||\varphi||_A$ for $\varphi\in A$ and $
|\varphi|\leq 1.$ 

Proof: We have shown that $F(\psi)\geq 0$ for all $
\psi\in A$ and $\psi>0.$ It follows that for any two $\varphi, \varphi_1 \in
A$ with $\varphi< \varphi_1,$ $F(\varphi)\leq F(\varphi_1).$ In particular,
if we denote $k=F(1),$ then (choose $\varphi_1=0$ and $1,$ respectively) for
any $0<\varphi<1,$ $0\leq F(\varphi)\leq k.$ If $-1<\varphi<1$ with $
\varphi=0$ somewhere (i.e. $|\varphi|<1$), then $0<1+\varphi<2,$ and $0\leq
F(1+\varphi)\leq 2k.$ That is, $|F(\varphi)|\leq k.$ For $\varphi\in A,
|\varphi|\leq 1,$ we consider $\varphi/2.$ Then $|\varphi/2|<1$ and $
|F(\varphi/2)|\leq k,$ i.e., $|F(\varphi)|\leq 2k$ for $\varphi\in A$ and $
|\varphi|\leq 1.$ Therefore, $F$ is a (complex-valued) continuous linear
functional on $A.$ The rest of the proof is exactly the same as that in the first proof.

\section{Hausdorff-Young Inequality}

\begin{thm}[Riesz-Thorin Convexity Theorem]
 Assume that $p_0\neq p_1$, $q_0\neq q_1$ and let $T$ be a linear operator such that 
\[T: L_{p_0}(U,d\mu)\rightarrow L_{q_0}(V,d\nu) \]
with norm $M_0$, and that 
\[
T: L_{p_1}(U,d\mu)\rightarrow L_{q_1}(V,d\nu) 
\]
with norm $M_1.$ Then 
\[
T: L_{p}(U,d\mu)\rightarrow L_{q}(V,d\nu) 
\]
with norm $M_\theta\leq M_0^{1-\theta}M_1^{\theta},$ provided that $0\leq
\theta\leq 1$ and 
\[
\frac{1}{p}=\frac{1-\theta}{p_0}+\frac{\theta}{p_1}; \,\,\,\frac{1}{q}=\frac{
1-\theta}{q_0}+\frac{\theta}{q_1}. 
\]
\end{thm}

{\bf Proof:} (See [4])

\begin{thm}[Hausdorff-Young]
For any $f\in L_p(T)$, $1\leq p\leq 2,$ 
\[
||\hat{f}||_q\leq ||f||_p, 
\]
where $q$ is the exponent conjugate to $p$ and, of course, $||\hat{f}||_q$
is the norm of the sequence of Fourier coefficients of $f$ in $l_q.$
\end{thm}

{\bf Proof:}\quad The Hausdorff-Young inequality is a simple consequence of
the Riesz-Thorin convexity theorem.

Let $Tf=\hat{f}.$ Note that $Tf$ is the sequence of Fourier coefficients of $
f$. It is trivial that 
\[
||Tf||_{\infty}\leq ||f||_1. 
\]
In addition, Bessel's inequality gives 
\[
||Tf||_2\leq ||f||_2. 
\]
Given $p$ with $1\leq p\leq 2,$ let $0\leq \theta\leq 1$ be such that 
\[
\frac{1}{p}=\frac{1-\theta}{1}+\frac{\theta}{2}.
\]
By the Riesz-Thorin theorem, we have 
\[
||Tf||_q\leq ||f||_p,\,\,\forall f\in L_p, 
\]
where $q$ is given by 
\[
\frac{1}{q}=1-\frac{1}{p}. 
\]

It is worth noting that we couldn't get the best constant in the Hausdorff-
Young inequality. Beckner proved (Annals of Math, 102(1975)) for the Fourier
transforms on $R$ that 
\[
||\hat{f}||_q\leq \sqrt{\frac{p^{1/p}}{q^{1/q}}} ||f||_p. 
\]

The following proof of the Hausdorff-Young inequality is due to A.P.Calderon
and A. Zygmund. It suffices to show that for any trigonometric polynomial $f$
with Fourier coefficients $c=(c_n)$ and $||f||_p=1$ we have $||c||_q\leq 1.$
Using the duality, we see that it suffices to show that 
\[
|\sum c_nd_n|\leq 1 
\]
for every sequence $d$ with $||d||_p=1.$

Put $f(t)=F(t)^{1/p} E(t)$ for $t\in T$ such that $F(t)=|f(t)|^p\geq 0$ and $
|E(t)|=1.$ ($E(t)=exp\{i arg (f(t))\}$. In case $f(t)=0$, simply define $
E(t)=1$). Similarly, put $d_n=D_n^{1/p}e_n$ with $D_n\geq 0$ and $|e_n|=1.$

Using these functions we write $\sum c_nd_n$ as 
\[
\sum c_nd_n=\sum D_n^{1/p}e_n\int F(t)^{1/p}E(t)e^{-int}dt. 
\]

Introducing the complex variable $z$, we define the function 
\[
Q(z)=\sum D_n^{z}e_n\int F(t)^{z}E(t)e^{-int}dt. 
\]

Using the Lebesgue Dominated Convergence Theorem we can prove that $Q(z)$ is
analytic in $Re z>0.$ 

Since the sum has only finitely many terms, each one (as function of $z$) is
bounded in the strip $\frac{1}{2}\leq Re z\leq 1.$ Hence $Q(z)$ is bounded
in this strip with bound depending on $d_n^{\prime}s$ and $f$.

For $Re z=1,$ we have 
\[
|Q(1+it)|\leq \sum D_n\int F(t) dt=1. 
\]
For $Re z=\frac{1}{2},$ the Schwarz inequality gives 
\[
|Q(\frac{1}{2}+i\theta)|\leq (\sum D_n)^{1/2} (\sum |\int F(t)^{\frac{1}{2}
+i\theta}E(t) e^{-int} dt|^2)^{1/2}. 
\]
The integral is the Fourier coefficient of $F(t)^{\frac{1}{2}+i\theta}E.$
Bessel's inequality gives that 
\[
(\sum |\int F(t)^{\frac{1}{2}+i\theta}E(t) e^{-int} dt|^2)^{1/2}\leq
||f^{p/2}||_2=(||f||_p)^{p/2}=1. 
\]
Therefore 
\[
|Q(\frac{1}{2}+i\theta)|\leq 1. 
\]

$Q(z)$ is analytic and bounded in the strip $\frac{1}{2}\leq Re z\leq 1$
(with bound depending on $d_n^{\prime}s$ and $f$) and bounded by $1$ on the
lines $Re z=\frac{1}{2}$ and $Re z=1.$ By Hadamard's three-lines theorem, $
|Q(z)|\leq 1$ for all $z$ throughout the strip. In particular, taking $z=
\frac{1}{p}$ in $Q(z),$ we have $|\sum c_nd_n|\leq 1$ for every sequence $d$
with $||d||_p=1.$ 

\vspace{0.5in}

\begin{thm}
Let $f$ be a summable function whose coefficient sequence is in $l^p$, $
1<p<2.$ Show that $f\in L^q$ and $||f||_{q}\leq ||\hat{f}||_p,$ where
$\frac{1}{p}+\frac{1}{q}=1.$
\end{thm}

{\bf Proof:}\quad Let $1<p<2$ and let $c=(c_k)_{k=-\infty}^{\infty}\in l^p.$
For $t\in (-\pi, \pi]$, we define 
\[
(T^*c)(t)=\sum c_ke^{ikt}. 
\]
We must show that $T^*$ is well-defined in the sense that 
\[
s_n(t)=\sum_{k=-n}^nc_ke^{ikt} 
\]
converges to a function $f(t)$ on $(-\pi,\pi]$ in the norm of $L^q,$ where $
q $ is the exponent conjugate of $p.$

For any $h\in L^p,$ by the H\"{o}lder and Hausdorff-Young inequalities, 
\begin{eqnarray*}
|\int_{-\pi}^{\pi}h(t)\overline{s_n(t)}d\sigma| &=&|\sum_{k=-n}^n\hat{h}(k)
\overline{c_k}| \\
&\leq & (\sum_{k=-n}^n|\hat{h}(k)|^q)^{1/q}(\sum_{k=-n}^n|c_k|^p)^{1/p} \\
&\leq &||\overline{h}||_q ||c||_p\leq ||h||_p||c||_p.
\end{eqnarray*}
This implies that $||s_n||_q\leq ||c||_p$ for all $n.$ Note that this is valid
for any $c\in l^p.$ We have 
\[
||s_m-s_n||_q=||\sum_{n<|k|\leq m} c_ke^{ikx}||_q\leq \sum_{n<|k|\leq m}
|c_k|^p. 
\]
Therefore, $s_n$ is a Cauchy sequence in $L^q$ and hence there exists an $
f\in L^q$ so that $||s_n-f||_q\rightarrow 0.$ We simply define $
(T^*c)(t)=f(t).$ Note that $T^*$ is an adjoint operator to $T$, the finite
Fourier transform, in the sense that 
\[
<T(h), c>=<h, T^*(c)>, 
\]
for all $h\in L^p$ and $c\in l^p,$ where $<T(h), c>=\sum \hat{h}(k)\overline{
c_k}$ and $<h, T^*(c)>=\displaystyle{\int h(t)\overline{f}(t)dt}$ with $f$
defined as the $L^q$ limit of $s_n.$

Moreover, for each $k$ and for any $n>|k|$, by H\"{o}lder's inequality, 
\[
|\hat{f}(k)-c_k|=|\int (f(t)-s_n(t))e^{-ikt} d\sigma|\leq ||f-s_n||_q. 
\]
Therefore, $\hat{f}(k)=c_k.$ 

{\bf Remarks:}

\begin{enumerate}
\item The case $p=2$ is the theorem of Riesz-Fischer.

\item The case $p=1$, to every $c\in l^1$ we may assign the continuous
function $f(t)=\sum c_ke^{ikt}.$ Since the series converges uniformly, $c_k=
\hat{f}(k)$ and $||f||_C\leq ||c||_1.$

\item The restriction of the theorem to $1\leq p\leq 2$ is essential. For
there is a sequence $c\in l^q$ for all $q>2$ and yet is not the finite
Fourier transform of any function in $L^1.$

\vspace{0.1in}

The series
\[
\sum \pm n^{-1/2}\cos nx, 
\]
with a suitable choice of signs, is a desired example as shown by the
following theorem: If $\sum (a_n^2+b_n^2)$ diverges, then almost
all the series 
\[
\sum r_n(t)(a_n\cos nx+b_n\sin nx) 
\]
are not Fourier series (because almost all the series are almost everywhere
non-Fej\'{e}r summable).
\end{enumerate}

\begin{thm}
The restriction of the Hausdorff-Young inequality to $1\leq p\leq 2$ is
essential, for there is a continuous function $f\in C$ (hence $f\in L^p$ for
all $p>0$) such that $||\hat{f}||_q=\infty$ for all $q<2.$ Therefore, it is
impossible that for some $p>2$, we would have $||\hat{f}||_q\leq ||f||_p$
for $f\in L_p.$
\end{thm}

{\bf Proof:}\quad The construction of the desired function follows from the
following theorem (see [5]): If $\sum a_n^2+b_n^2<\infty,$ then
almost all series 
\[
\sum r_n(t)(a_n\cos nx+b_n\sin nx) 
\]
converges almost everywhere in $x\in [0,2\pi].$ If $\sum (a_n^2+b_n^2)(log
n)^{1+\epsilon}<\infty$ for some $\epsilon>0,$ then almost all series 
\[
\sum r_n(t)(a_n\cos nx+b_n\sin nx) 
\]
converges uniformly and so are Fourier series of continuous functions.

The series 
\[
\sum \pm\frac{\cos nx}{n^{1/2}\log ^2n} 
\]
is, for a suitable choice of signs, a case in point.

\begin{thm}[Hadamard's three-lines theorem]
Assume that $f(z)$ is analytic on the open strip $0<Re z<1$ and bounded (by $
B$) and continuous on the closed strip $0\leq Re z\leq 1.$ If $|f(it)|\leq
M_0$ and $|f(1+it)|\leq M_1,$ $-\infty<t<\infty,$ then we have $|f(\theta
+it)|\leq M_0^{1-\theta}M_1^{\theta}.$
\end{thm}

{\bf Proof:}\quad We assume first that $M_0=M_1=1.$ We have to prove that $
|f(z)|\leq 1$ for all $z$ in the strip.

For each $\epsilon>0$, we define 
\[
h_\epsilon(z)=\frac{1}{1+\epsilon z}, \quad z\in \mbox{the strip}. 
\]
Since $Re (1+\epsilon z)\geq 1$ in the closed strip, we have $|h_\epsilon|<1$
in the closed strip, so that 
\[
|f(z)h_\epsilon(z)|\leq 1 
\]
for $z$ in the boundaries of the strip. Also $|1+\epsilon z|\geq \epsilon
|y|,$ so that 
\[
|f(z)h_\epsilon (z)|\leq \frac{B}{\epsilon |y|},\quad z\in 
\mbox{the closed
strip}.
\]

Let $R$ be the rectangle cut off from the closed strip by the lines $y=\pm
B/\epsilon.$ Since $|fh_\epsilon|\leq 1$ on the boundaries of $R,$  $
|fh_\epsilon|\leq 1$ on $R,$ by the maximum modulus principle. But the above also shows
that $|fh_\epsilon|\leq 1$ on the rest of the closed strip. Thus $
|fh_\epsilon|\leq 1$ throughout the closed strip. If we fix $z$ in the strip
and let $\epsilon\rightarrow 0,$ we obtain $|f(z)|\leq 1.$

We now turn to the general case. Put 
\[
g(z)=M_0^{1-z}M_1^z, 
\]
where $M_i^{\zeta}=exp\{\zeta\log M_i\}$ for complex $\zeta.$ Then $g(z)$ is
entire, $g$ has no zero, $1/g$ is bounded in the closed strip, 
\[
|g(it)|=M_0,\qquad |g(1+it)|=M_1, 
\]
and hence $f/g$ satisfies our previous assumptions. Thus $|f/g|\leq 1$ in
the strip, and this gives $|f(\theta +it)|\leq M_0^{1-\theta}M_1^{\theta}$
for all $0\leq\theta\leq 1.$ 

\vspace{0.1in}

{\bf An Alternative Proof:}

Let $\epsilon>0$ and $\lambda\in R.$ Define 
\[
F_\epsilon (z)=exp\{\epsilon z^2+\lambda z\}F(z). 
\]
Then 
\[
F_\epsilon(z)\rightarrow 0, \,\, \mbox{as}\,\, Im z\rightarrow \pm \infty 
\]
and 
\[
|F_\epsilon(it)|\leq M_0,\,\, |F_\epsilon(1+it)|\leq M_1
e^{\epsilon+\lambda}. 
\]
By the Phragmen-Lindel\"{o}f principle we therefore obtain 
\[
|F_\epsilon(z)|\leq max\{M_0, M_1e^{\epsilon+\lambda}\}.
\]
That is, 
\[
|F(\theta+it)|\leq exp\{-(\theta^2-t^2)\} max\{M_0e^{-\theta\lambda},
M_1e^{(1-\theta)\lambda+\epsilon}\}. 
\]
This holds for any fixed $\theta$ and $t.$ Letting $\epsilon\rightarrow 0$
we conclude that, if $\rho=exp\{\lambda\},$ 
\[
|F(\theta+it)|\leq max\{M_0\rho^{-\theta}, M_1\rho^{1-\theta}\}. 
\]

\section{A Theorem of Minkowski}

Let 
\[
T^2=\{(e^{i2\pi x}, e^{i2\pi y}): x, y \in R\}. 
\]
$T^2$ is called the 2-dimensional torus, which is the Cartesian product of
the unit circle $T=\{e^{i2\pi x}: x\in R\}.$

Let $(m,n)$ be a lattice (integer coordinates) point in the plane and let $
f(x,y)$ be a summable function on the unit square 
\[
E=\{(x,y): 0<x<1, 0<y<1\} 
\]
with extension to $R^2$ periodically. Define the Fourier coefficients 
\[
a_{m,n}(f)=\int\int_Ef(x,y)e^{2\pi i(mx+ny)}dxdy. 
\]

We may prove the Parseval relation on $L^2(T^2).$ For any trigonometric
polynomial 
\[
p(x,y)=\sum\sum_{\mbox{finite sum}} a_{m,n}e^{i2\pi (mx+ny)}, 
\]
the Parseval relation holds. Since such trigonometric polynomials are dense
in $L^2(T^2)$, the Parseval relation holds for every $f\in L^2(T^2).$

\begin{thm}[Minkowski]
Let $C$ be a convex body in $R^d$ of volume $V$ and symmetric about the
origin. If $V>2^d,$ then $C$ contains a lattice point other than the origin.
\end{thm}

{\bf Proof:}\quad We work with $d=2.$ Let $C$ be a convex body in $R^2$ of
volume $V$ and be symmetric about the origin. Assume that $C$ contains no 
lattice point except the origin. We want to show that the area $V$ of $C$ is 
$\leq 4.$

Let $\phi(x,y)$ be the characteristic function of $C.$ Let 
\[
f(x,y)=\sum_{m,n}\phi(2(x-m),2(y-n)). 
\]
Assume that $C$ is bounded (If $C$ is not bounded, we consider the
intersection of $C$ and the circle with center at the origin and radius $R.$
If $V>4,$ then for some large enough $R$ the area of the intersection is $>$
4. Also note that the intersection is convex, symmetric, and bounded. Thus
the Minkowski theorem applies). For each $(x,y),$ this sum has only finitely
many nonzero terms.

Let $E$ be the unit square defined as above. The Parseval relation asserts
that 
\[
\sum_{m,n}|a_{m,n}(f)|^2=\int\int_E |f(x,y)|^2 dx dy. 
\]

We calculate $a_{m,n}$ as follows: 
\begin{eqnarray*}
a_{m,n}&=&\int\int_E f(x,y)e^{-2\pi i(mx+ny)} dx dy \\
&=&\int\int_E \sum_{m,n}\phi (2(x-m), 2(y-n))e^{-2\pi i(mx+ny)} dx dy \\
&=&\int\int_{R^2}\phi (2x, 2y)e^{-2\pi i(mx+ny)} dx dy \\
&=&2^{-d}\int\int_{C} \phi (x, y)e^{-\pi i(mx+ny)} dx dy.
\end{eqnarray*}

The last equality is simply the result of change of variables. For the one above the last
equality, we denote by $E_{-m,-n}$ the square with the lower left corner $
(-m,-n).$ Then we have: 
\begin{eqnarray*}
&\empty &\int\int_{R^2} \phi (2x, 2y)e^{-2\pi i(mx+ny)} dx dy \\
&=&\sum_{m,n}\int\int_{E_{-m,-n}}\phi (2x, 2y)e^{-2\pi i(mx+ny)} dx dy \\
&=&\sum_{m,n}\int\int_{E} \phi (2(x-m), 2(y-n))e^{-2\pi i(mx+ny)} dx dy \\
&=&\int\int_{E} \sum_{m,n}\phi (2(x-m), 2(y-n))e^{-2\pi i(mx+ny)} dx dy.
\end{eqnarray*}

On the other hand, we calculate $\int\int_E|f(x,y)|^2 dx dy.$ 
\begin{eqnarray*}
\int\int_E|f(x,y)|^2 dx dy&=& \int\int_E f(x,y)\sum_{m,n}\phi (2(x-m),
2(y-n))dx dy \\
&=&\int\int_{R^2} f(x,y)\phi (2x, 2y) dx dy \\
&=&\int\int_{R^2} \sum_{m,n}\phi(2(x-m), 2(y-n))\phi (2x, 2y) dx dy \\
&=&\sum_{m,n}\int\int_{R^2} \phi(2(x-m), 2(y-n))\phi (2x, 2y) dx dy \\
&=&2^{-d}\sum_{m,n}\int\int_{R^2} \phi(x-2m, y-2n)\phi (x, y) dx dy \\
&=&2^{-d}\sum_{m,n}\int\int_{C} \phi(x-2m, y-2n) dx dy.
\end{eqnarray*}

The Parseval relation gives rise to 
\[
2^{-2d}|\int\int_C\phi(x,y)e^{-\pi i(mx+ny)}dx dy|^2=
=2^{-d}\sum_{m,n}\int\int_{C} \phi(x-2m, y-2n) dx dy. 
\]
If $C$ contains no lattice point except the origin, then one can show that
for $(x,y)\in C$ and $(m,n)\neq (0,0),$ $(x-2m,y-2n)\not\in C$ and so every
term in the sum is zero except the one with $m=n=0.$  The term with $(m,n)=(0,0)$ equals $2^{-d}V.$ Thus we
have 
\[
2^{-d}\sum_{m,n}|\int\int_C\phi(x,y)e^{-\pi i(mx+ny)}dx dy|^2=V. 
\]
The term on the left with $(m,n)=(0,0)$ is $2^{-d}V^2,$ and therefore, $
2^{-d}V^2\leq V$, that is, $V\leq 2^d.$

\vspace{0.1in}

\begin{thm}
If $C$ is a convex body in $R^d$ of volume $V=2^d$ and symmetric about the
origin, then there is a lattice point $(m,n)\neq (0,0)$ in $C$ or on its
boundary.
\end{thm}

{\bf Proof:}\quad Assume that, by a contradiction, $\overline{C}$ contains
no lattice point other than the origin. We assume that $C$ is bounded. Thus $
\overline{C}$ is compact and there is $\delta>0$ such that $d(p,
C)\geq\delta>0$ for all lattice points $p$ other than the origin. We may
expand $\overline{C}$ slightly to a subset $D$ of $R^d$ so that $D$ is
convex and symmetric about origin and yet contains no lattice point other
than the origin. Since the volume of $D$ is $>2^d,$ this is in contradiction
to Minkowski's theorem.

It remains to show the construction of $D.$ Let 
\[
D=\{x\in R^d: dist (x,C)\leq\frac{\delta}{2}\}. 
\]
We claim that if $C$ is (closed) convex, then so is $D$. Let $x$ and $y\in D$
(Assume that they are not in $C.$ Otherwise, nothing needs to be done.) Let $x_0$
and $y_0\in C$ such that $|x-x_0|=dist (x,C)$ and $|y-y_0|=dist (y,C).$ For $
0\leq \lambda\leq 1,$ we have $|(\lambda x+(1-\lambda)y )- (|\lambda
x_0+(1-\lambda)y_0)| \leq \lambda |x-x_0|+(1-\lambda )|y-y_0|\leq \frac{
\delta}{2}.$ Note that $C$ is convex, $\lambda x_0+(1-\lambda)y_0\in C.$
Therefore, $dist (\lambda x+(1-\lambda)y, C)\leq \frac{\delta}{2}$ and so $
\lambda x+(1-\lambda)y\in D.$ Clearly, if $C$ is symmetric about the origin
then so is $D.$ 

\vspace{0.1in}

\begin{thm}
Prove that if $k>0$ and $a,b,c,$ $d$ are real with $|ad-bc|\leq 1,$ then
there are integer pairs $(m,n)\neq (0,0)$ such that 
\[
|am+bn|\leq k\qquad |cm+dn|\leq k^{-1}. 
\]
Deduce that for every real number $a,$ there are infinitely many integer
pairs $(m,n)$ such that 
\[
|a+\frac{n}{m}|\leq \frac{1}{m^2}. 
\]
\end{thm}

{\bf Proof:}\quad Let 
\[
C=\{(x,y): |ax+by|\leq k \,\,\mbox{and}\,\, |cx+dy|\leq \frac{1}{k}\}. 
\]
Clearly, $C$ is the parallelogram centered at $(0,0)$ and bounded by $
ax+by+k=0,$ $ax+by-k=0,$ $cx+dy+\frac{1}{k}=0,$ and $cx+dy-\frac{1}{k}=0.$
Thus, $C$ is convex and symmetric. Note that $C$ is the set of all points
whose distances are $\leq \frac{k}{\sqrt{a^2+b^2}}$ from $ax+by=0$ and are $
\leq \frac{1}{k\sqrt{c^2+d^2}}$ from $cx+dy=0.$ To find the area of $C$ we
observe that the mapping $x=du-bv$ and $y=-cu+av$ maps $C$ to the rectangle
in the $u-v$ plane bounded by $(ad-bc)u\pm k=0$ and $(ad-bc)v\pm \frac{1}{k}
=0.$ The area of the rectangle is 
\[
\frac{2y}{|ad-bc|}\cdot \frac{2}{k|ad-bc|}. 
\]
The Jacobian of the mapping is $|ad-bc|.$ Thus 
\[
V=\frac{2y}{|ad-bc|}\cdot \frac{2}{k|ad-bc|}\cdot |ad-bc|= =\frac{4}{|ad-bc|}
. 
\]
If $|ad-bc|\leq 1$ then $V\geq 4$ and so there is a lattice point $(m,
n)\neq (0,0)$ in $C$ or on its boundary. Note that $C$ is closed. The
lattice point is in $C$ anyway.

\vspace{0.1in}

{\bf An alternative proof:} Let $a\in R.$ We want to show there are
infinitely many integer pairs $m$ and $n$ such that 
\[
|a-\frac{n}{m}|\leq \frac{1}{m^2}. 
\]

Let $N$ be a positive integer. We divide $[0,1]$ into $N+1$ parts so that each part has length $
\frac{1}{N}.$ Consider $\{ja\}$, the fraction part of $ja$, for $
j=0,1,\cdots, N.$ By pigeonhole principle, there are $j$ and $k$ such that
both $\{ja\}$ and $\{ka\}$ lie in some interval of length $\frac{1}{N}$. It
follows that there is an integer $n$ so that $|ja-ka-n|\leq\frac{1}{M}.$ Let 
$|j-a|=m$. Then 
\[
|a-\frac{n}{m}|\leq \frac{1}{mN}\leq \frac{1}{m^2}. 
\]

\section{Measures with bounded powers}

Let $M(R)$ be the algebra of complex bounded Borel measures on $R.$ The
multiplication of $\mu$ and $\nu\in M(R)$ is defined as the convolution $
\mu*\nu,$ which is the measure in $M(R)$ satisfying 
\[
\int h(t)d(\mu*\nu)(t)=\int\int h(s+t) d\mu(s)d\nu(t),\qquad \forall h\in
C_0(R). 
\]
The measure $\delta(t)$ with unit mass at $0$ is an identity for this
algebra. If $\mu$ is the measure with mass $\epsilon$, $|\epsilon|=1,$ at $x,$
then $\nu$, defined as the measure with mass $\overline{\epsilon}$ at $-x$,
is an inverse of $\mu.$ To see this, let $h\in C_0(R),$ we have: 
\begin{eqnarray*}
\int h(t)d(\mu*\nu)(t)&=&\int\int h(s+t) d\mu(s)d\nu(t) =\int h(x+t)\epsilon
d\nu(t) \\
&=&h(0)\epsilon \overline{\epsilon}=h(0) =\int h(t)d\delta(t).
\end{eqnarray*}

If $\mu$ is a point measure with mass $\epsilon$ at $x$ then the power $
\mu^{*n}$ (in the sense of convolution) is the point mass with measure $
\epsilon^n$ at $nx.$ To see this, let $\nu=\mu^{*(n-1)}$ and $h\in C_0(R)$.
We then have 
\begin{eqnarray*}
\int h(t)d(\mu*\nu)(t)&=&\int\int h(s+t) d\mu(s)d\nu(t) =\int h(x+t)\epsilon
d\nu(t) \\
&=&h(x+(n-1)x)\epsilon \epsilon^{n-1} =h(nx)\epsilon^n=\int h(t)d\mu^{*n}(t).
\end{eqnarray*}
Also note that if $n$ is negative, then $\mu^{*n}$ is defined as 
\[
\mu^{*n}=(\mu^{-1})^{*|n|}. 
\]
Of course, in order for this definition to make sense, we must agree that
whenever we mention $\mu^{*n}$ with negative $n$, we admit that $\mu$ has
inverse.

We can identify the set of all bounded point-measures with mass at integers
with $l^1$ by considering such measures as functions on $Z.$ With this
correspondence, the identity $\delta$ corresponds to $\delta
=\{\delta(n)\}_{n=-\infty}^\infty \in l^1$ with $\delta(0)=1,$ and $
\delta(n)=0$ for $n\neq 0,$ and the convolution of $\mu$ and $\nu$ is simply
the convolution of two elements $\mu$ and $\nu$ in $l^1$ defined as 
\[
(\mu*\nu)(j)=\sum_{k=-\infty}^\infty \mu(k)\nu(j-k). 
\]

When does a given element $\mu\in l^1$ have an inverse? That is, given $
\mu\in l^1,$ does there exist $\nu\in l^1$ such that $\mu*\nu=\delta?$ To
answer this question we need Wiener's theorem. By Weiner's theorem we conclude that if $\mu\in l^1$
and $m(x)=\sum \mu(j)e^{ijx}$ is nowhere zero, then there is $\nu\in l^1$
such that $\mu*\nu=\delta.$ In fact, the Fourier transform of ${\displaystyle
\frac{1}{m(x)}}$ will be the desired $\nu.$

The sequence $a=(a_n)\in l^1$ with $a_0=1$ and $a_1=1$ has no inverse. In
fact, if there were $b\in l^1$ such that $a*b=\delta,$ then we would have $
\mathcal{F}(a*b)= \mathcal{F}(a)\cdot \mathcal{F}(b)=1$ with respect to the ordinary
multiplication. This is a contradiction, because the two functions are
continuous and $\mathcal{F}(a)=1+e^{ix}$ equals zero at $x=\pi$, $\mathcal{F}
(a)\cdot \mathcal{F}(b)\neq 1$ at $x=\pi$ for any value of $\mathcal{F}(b)$ at $
x=\pi.$ This shows that if $\mu$ has an inverse, then $\mathcal{F}(\mu)$ can be
zero nowhere on $T.$

Moreover, with this correspondence, if $\mu$ is the measure with point mass $
\epsilon$ at $p$ ($p$ is an integer) i.e., $\mu(p)=\epsilon$ and $\mu(k)=0$
for $k\neq p,$ then $\mu^{*n}$ is the measure with point mass $\epsilon^n$
at $np.$ To see this, let $\nu=\mu^{*(n-1)},$ then 
\[
\mu^{*n}(k)=\sum_{l=-\infty}^\infty \mu(l)\nu(k-l)=\mu(p)\nu(k-p) 
\]
$=0$ for all $k\neq np;$ $\epsilon^n$ for $k=np.$

Let $\mu$ be a point measure with mass at integers and let 
\[
m(x)=\sum \mu(k)e^{ikx}. 
\]
$m(x)$ is the inverse Fourier transform of $\mu(n),$ and the Fourier
coefficient $\displaystyle{\int m(x)e^{-ikx}d\sigma (x)}$ of $m(x)$ is $
\mu(k),$ for $k=0,\pm 1, \pm 2, \cdots.$

One can show that 
\[
m^{n}(x)=\sum \mu^{*n}(k)e^{ikx}, 
\]
where $m^n(x)$ is the nth power of $m$ with respect to ordinary
multiplication. We prove this with $n=2.$ Using the Cauchy product of two
series, we have 
\begin{eqnarray*}
\int m^2(x)e^{-ikx}d\sigma (x)&=& \int (\sum \mu(m)e^{imx})(\sum
\mu(n)e^{inx})e^{-ikx}d\sigma (x) \\
&=& \int \biggl(\sum_n \sum_{l=-\infty}^{\infty}\mu(l)\mu(n-l)e^{inx}\biggr) 
e^{-ikx}d\sigma (x) \\
&=& \sum_{l=-\infty}^{\infty}\mu(l)\mu(k-l).
\end{eqnarray*}
Therefore, 
\[
m^2(x)=\sum_k \biggl(\sum_{l=-\infty}^{\infty}\mu(l)\mu(k-l)\biggr)e^{ikx}
=\sum_k\mu^{*2}(k)e^{ikx}. 
\]

\begin{lem}
If $f\in l^1$ with $||f||_1\leq K$ and $||f||_2=1,$ then $||f||_4\geq r,$
where $r=K^{-1/2}>0.$
\end{lem}

{\bf Proof:}\quad For simplicity of notation, we view $f$ as a function
defined on $R$ equipped with the measure $\mu$ having unit mass at each
integer. Let $0<\theta<1$ be such that 
\[
\frac{1}{2}=\frac{1-\theta}{1}+\frac{\theta}{4}. 
\]
In fact, $\theta=\frac{2}{3}$ will do. Using H\"{o}lder's inequality with
indices $p={\displaystyle\frac{1}{2(1-\theta)}}$ and $q={\displaystyle\frac{2
}{\theta}},$ we get 
\[
(\int_R |f|^2d\mu)^{1/2}= (\int_R
|f|^{2(1-\theta)}|f|^{2\theta}d\mu)^{1/2}\leq (\int_R |f|
d\mu)^{1-\theta}(\int |f|^4d\mu)^{\theta/4}. 
\]
The desired estimation for $||f||_4$  follows.

\vspace{0.3in}

{\bf An alternative proof:}

\vspace{0.1in}

Using Schwarz' inequality twice, we get 
\begin{eqnarray*}
1&=&\|\mu\|_2^2\leq \||\mu|^{3/2}\|_2 \||\mu|^{1/2}\|_2 \\
&=& \||\mu|^2|\mu|\|_1^{1/2} \|\mu\|_1^{1/2} \\
&\leq& (\||\mu|^2\|_2^{1/2}\|\mu\|_2^{1/2})^{1/2} K^{1/2} \\
&=&\||\mu|\|_4 K^{1/2}.
\end{eqnarray*}

\begin{thm}[Beurling-Helson]
$\mu\in l^1$ has bounded powers, i.e., $||\mu^{*n}||_{l^1}\leq K$ for all $
n=0,\pm 1,\cdots,$ if and only if it satisfies $|\mu(p)|=1$ for some integer 
$p$, $\mu(m)=0$ for all $m\neq p.$
\end{thm}

{\bf Proof:}\quad (1). Prove that $m(x)=\mathcal{F}^{-1}(\mu)$ can be
written as $e^{i\phi(x)}$ for some $\phi(x).$ Reduce the theorem to proving $
\phi(x)=px$ for some integer $p$.

Let $m(x)=\sum \mu(k)e^{ikx},$ $x\in [0,2\pi).$ Then 
\[
m^n(x)=\sum \mu^{*n}(k)e^{ikx}. 
\]
If $||\mu^{*n}||_{l^1}\leq K$ for all $n=0,\pm 1,\cdots,$ then $|m^n(x)|\leq
K$ for all $n$ and all $x.$ Therefore, $|m(x)|=1$ for all $x.$ Since $m(x)$
is continuous, we can write 
\[
m(x)=e^{i\phi(x)}, 
\]
where $\phi(x)$ is continuous and $\phi(x+2\pi)-\phi(x)$ is a (fixed
constant) multiple of $2\pi$ for all $x.$
 To emphasize, we write again that 
\[
m^n(x)=e^{in\phi(x)}=\sum \mu^{*n}(k)e^{ikx}.
\]
The conclusion of the theorem is equivalent to $\phi(x)=px+b$ for some
integer $p$ and real number $b.$ Without loss of generality, we assume that 
$\phi(0)=0,$ and we shall prove $\phi(x)=px.$

(2).  Let $\Phi(r,s,t)=\phi(t-r)+\phi(r)-\phi(t-s)-\phi(s).$ Show that $
e^{i\Phi(r,s,t)}=w$ for some constant $w$ on a set of positive measure in $
[0,2\pi)^3.$ (If $\phi(x)=px,$ then $\Phi(r,s,t)=0$
identically.)

Note that ${\displaystyle\sum |\mu^{*n}(k)|\leq K}$ (by assumption) and
that, by the Parseval relation, 
\[{\displaystyle\sum_k |\mu^{*n}(k)|^2=\int
|m^n(x)|^2 d\sigma(x)=1}.\] Therefore, by the lemma, ${\displaystyle \sum
|\mu^{*n}(k)|^4\geq r^4>0}.$

On the other hand, since 
\[
\mathcal{F}(\int e^{in(\phi(t-s)+\phi(s))}d\sigma(s)) = \biggl\{(\mu^{*n}(k))^2
\biggr\}_{k\in Z}, 
\]
by the Parseval relation we have 
\[
\sum |\mu^{*n}(k)|^4= \int \biggl|\int e^{in(\phi(t-s)+\phi(s))}d\sigma(s) 
\biggr|^2d\sigma(t). 
\]
The integral can be written as 
\begin{eqnarray*}
&.&\int \biggl|\int e^{in(\phi(t-s)+\phi(s))}d\sigma(s)\biggr|^2d\sigma(t) \\
&=&\int \biggl(\int e^{in(\phi(t-s)+\phi(s))}d\sigma(s)\biggr) \biggl(\int
e^{-in(\phi(t-r)+\phi(r))}d\sigma(r)\biggr)d\sigma(t) \\
&=& \int\int\int e^{in\Phi(r,s,t)}d\sigma(r,s,t),
\end{eqnarray*}
where $\Phi(r,s,t)=\phi(t-r)+\phi(r)-\phi(t-s)-\phi(s)$ and $
d\sigma(r,s,t)=d\sigma(r)d\sigma(s)d\sigma(t).$ Summarizing, we have 
\[
\int\int\int e^{in\Phi(r,s,t)}d\sigma(r,s,t) \geq r^4>0,\quad\forall n.
\qquad (1)\]

Using (1), we proceed to show that $e^{i\Phi(r,s,t)}=w$ for some value $w$
on a set of positive measure. Let $E\subset [0,2\pi)$ and Define 
\[
\nu(E)=\lambda(\{(r,s,t)\in [0,2\pi)^3 : 0\leq r,s,t <2\pi, \Phi(r,s,t)\in
E\}, 
\]
where $\lambda$ is the normalized Lebesgue measure on $[0,2\pi)^3.$ Then
\[
\int\int\int_{[0,2\pi)^3}
e^{in\Phi(r,s,t)}d\sigma(r,s,t)=\int_0^{2\pi}e^{inu}d\nu(u).\qquad (2) 
\]
Suppose that, by a contradiction, $\Phi(r,s,t)$ assumes each value only on a
null set (i.e., there is no subset of $[0,2\pi)^3$ of positive ($\lambda$
) measure on which $\Phi$ takes a constant value). Then, given any point $
\theta\in [0,2\pi),$ $\nu(\theta)=0$ and so $\nu$ has no point mass. By
Wiener's theorem 
\[
\lim_{N\rightarrow \infty}\frac{1}{2N+1}\sum_{k=-N}^N|\hat{\nu}(k)|^2 =0. 
\]
It follows that $\hat{\nu}(k)$ are not bounded away from zero and so that
there is a subsequence of $\hat{\nu}(k)$ converging to zero. Since $\hat{\nu}
(n)$ is given by the integral on the left hand of (2), by (1) we have $\hat{
\nu}(n)\geq r^4$ for all $n,$ and so any subsequence of $\hat{\nu}(n)$
cannot converge to zero. This is a contradiction. Therefore, there is a set
of positive measure, denoted by $A,$ such that $e^{i\Phi(r,s,t)}=w$ for some
value $w$ on $A.$

(3). Denote 
\[
A=\{(r,s,t)\in T^3: e^{i\Phi(r,s,t)}=w\}. 
\]
We will show that $A=T^3.$

 Let 
\[
\Psi(r,s,t)=\frac{w^{-1}e^{i\Phi(r,s,t)}+1}{2}. 
\]
Then $\Psi=1$ on $A$ and $|\Psi|<1$ anywhere else. ($w^{-1}e^{i\Phi}$ is a
complex number of modulus $1$ and so it is a point on the unit circle and ${
\displaystyle\frac{w^{-1}e^{i\Phi}+1}{2}}$ is a point inside the unit disk
unless $w^{-1}e^{i\Phi}=1$). It follows that 
\[
\Psi^n(r,s,t)\rightarrow \chi_A(r,s,t) 
\]
pointwise as $n\rightarrow\infty.$

To prove that $A=T^3,$ we show that we are able to extract a subsequence of $
\Psi^n$ converging to a continuous function on $T^3$. Then, since $
\Psi^n(r,s,t)\rightarrow \chi_A(r,s,t),$ we must have $A=T^3.$

Let $\psi\in l^1(Z^3)$ (sequence depending on three indices) be the Fourier
transform of $\Psi(r,s,t).$ Then $\psi^{*n}$ is the Fourier transform of $
\Psi^n$ (with respect to ordinary multiplication) and\,\,\,
 $||\psi^{*n}||_{l^1(Z^3)}\leq K^4$ for all $n.$

To see this, we first prove that 
\[
\|\mathcal{F}(e^{in\Phi(r,s,t)})\|_{l^1(Z^3)}\leq K^4,\qquad \forall n. 
\]
Note that 
\[
e^{in\Phi(r,s,t)}=e^{in\phi(t-r)}e^{in\phi(r)}
e^{-in\phi(t-s)}e^{-in\phi(s)}, 
\]
and that, by (1), 
\[
e^{in\phi(t-r)}=\sum_j \mu^{*n}(j)e^{-ji(t-r)}, 
\]
\[
e^{in\phi(r)}=\sum_{j} \mu^{*n}(j)e^{-jir}, 
\]
\[
e^{-in\phi(t-s)}=\sum_j \overline{\mu^{*n}(j)}e^{-ji(t-s)}, 
\]
and 
\[
e^{-in\phi(s)}=\sum_j \overline{\mu^{*n}(j)}e^{-jis}, 
\]
where each one is viewed as a Fourier series on $T^3.$ Therefore, the
Fourier series of $e^{in\Phi(r,s,t)}$ is a (Cauchy) product of all the
series above (or the convolution of all Fourier transforms) such that 
\begin{eqnarray*}
\|\mathcal{F}(e^{in\Phi(r,s,t)})\|&=& \|\mathcal{F}e^{in\phi(t-r)}* \mathcal{F}
e^{in\phi(r)}*\mathcal{F}e^{-in\phi(t-s)}*\mathcal{F}e^{in\phi(t-r)}\| \\
&\leq &\|\mathcal{F}e^{in\phi(t-r)}\|\cdot \|\mathcal{F}e^{in\phi(r)}\|\cdot\|
\mathcal{F}e^{-in\phi(t-s)}\| \cdot \|\mathcal{F}e^{-in\phi(s)}\|\leq K^4,
\end{eqnarray*}
where all the norms are taken in $l^1(Z^3).$

To see $||\psi^{*n}||_{l^1(Z^3)}\leq K^4$ for all $n,$ we write 
\begin{eqnarray*}
\Psi^n&=&\frac{1}{2^n}[(w^{-1}e^{i\Phi})^n
+c(n,1)(w^{-1}e^{i\Phi})^{n-1}+\cdots +1] \\
&=& \frac{1}{2^n}[w^{-n}e^{in\Phi}+c(n,1)w^{-(n-1)}e^{i(n-1)\Phi} +\cdots
+1].
\end{eqnarray*}
Therefore, 
\[
\mathcal{F}(\Psi^n)= \frac{1}{2^n}[w^{-n}\mathcal{F}(e^{in\Phi})+c(n,1)w^{-(n-1)}
\mathcal{F} (e^{i(n-1)\Phi}) +\cdots +1]. 
\]
It follows, since $|w|=1,$ that 
\begin{eqnarray*}
\|\psi^{*n}\|&=&\|\mathcal{F}(\Psi^n)\|_{l^1}\leq \frac{1}{2^n}[\|\mathcal{F}
(e^{in\Phi})\|+c(n,1)\|\mathcal{F} (e^{i(n-1)\Phi})\| +\cdots +1 \\
&\leq& K^4\frac{1}{2^n}[1+c(n,1)+\cdots+c(n,n-1)+1]=K^4.
\end{eqnarray*}

Since $(c_0(Z^3))^*=l^1(Z^3),$ by Alaoglu's theorem, $
\{\psi^{*n}\}$ is (sequentially) compact in the weak* topology of $l^1(Z^3),$ that is, there is $\rho\in l^1(Z^3)$ such that for some
subsequence $\psi^{*n_k},$ \, $\langle\psi^{*n_k}, c\rangle\rightarrow
\langle\rho, c\rangle$ for any $c\in c_0(Z^3)$ as $n_k\rightarrow\infty.$

Note that $|\Psi|\leq 1 $ so that $|\Psi^n|\leq 1$ and $\Psi^{n_k}
\rightarrow \chi_A$ pointwise as $k\rightarrow\infty.$ By the bounded
convergence theorem, we have 
\begin{eqnarray*}
\psi^{*n_k}(m,n,l)&=& \mathcal{F}(\Psi^{n_k})(m,n,l) \\
&=&\int\,\int\,\int_{T^3}\Psi^{n_k}(r,s,t) e^{i(mr+ns+lt)}drdsdt\rightarrow 
\mathcal{F}(\chi_A)(m,n,l),
\end{eqnarray*}
as $k\rightarrow\infty,$ that is, $\psi^{*n_k}\rightarrow \mathcal{F}(\chi_A)$
componentwise as $k\rightarrow \infty.$ On the other hand, $
\langle\psi^{*n_k}, c\rangle \rightarrow \langle\rho, c\rangle$ for any $
c\in c_0(Z^3)$ as $k\rightarrow\infty.$ Taking $c=e^{(m,n,l)}$, we have $
\psi^{*n_k}(m,n,l)\rightarrow \rho(m,n,l),$ i.e., $\psi^{*n_k}\rightarrow
\rho$ componentwise. Therefore, $\mathcal{F}(\chi_A)=\rho$ and so 
\[
\chi_A=\mathcal{F}^{-1}(\rho), \quad a.e. 
\]
Since the transform $\mathcal{F}^{-1}(\rho)$ of $\rho$ is continuous (because $
\rho\in l^1(Z^3)$), $A=T^3,$ that is, $e^{i\Phi(r,s,t)}$ assumes value $w$
everywhere on $T^3.$

\vspace{0.1in}

(4). Show that $\phi(t)=pt$ for some integer $p.$

Since $e^{i\Phi(r,s,t)}$ assumes value $w$ everywhere on $T^3,$ $
\Phi(r,s,t)=2k\pi +c,$ where $c=arg(w)$ and $k$ is an integer that might
depend on $r,s,t,$ apparently. Since $\Phi$ is continuous, $
\Phi(r,s,t)=2k\pi+c$ for all $r,s,t$ with a fixed constant $k.$ It follows
that $\phi(t-r)+\phi(r)$ is independent of $r.$ Since $\phi(0)=0,$ $
\phi(t-r)+\phi(r)=\phi(t).$ In other words, 
\[
\phi(t+s)=\phi(t)+\phi(s) 
\]
for all real $r,s.$ Therefore, $m=e^{i\phi}$ is a character of $R$  and, since $
\phi $ is periodic, $m$ is a character of $T.$ Thus, $\phi(t)=pt$ for some
integer $p.$ 

\vspace{0.2in}

{\bf Homomorphisms of $l^1$ into $l^1$}

\vspace{0.1in}

Let $h$ be a homomorphism of $l^1$ into itself, i.e., $h$ is a linear
mapping and $h(\mu*\rho)=h(\mu)*h(\rho)$ for all $\mu, \rho\in l^1.$ Note
that $h$ is continuous.

Let $e^{(n)},$ $n=0,\pm 1,\pm 2, \cdots,$ be a sequence of elements of $l^1$
such that $e^{(n)}(n)=1$ and $e^{(n)}(k)=0$ for $k\neq n.$ Note that $
e^{(0)} $ is the identity for the algebra $l^1.$ We have that $
h(\mu)=h(e^{(0)})*h(\mu)$ for all $\mu\in l^1.$ Therefore, $
h(e^{(0)})=e^{(0)}.$ Let 
\[
h(e^{(1)})=\mu\in l^1. 
\]
Note that $e^{(m)}*e^{(n)}=e^{(m+n)}.$ By the multiplicativity of $h$, $
h(e^{(n)})=\mu^{*n},$ and $\|\mu^{*n}\|_{l^1}\leq \|h\|$ (since $
\|e^{(n)}\|_{l^1}=1$) for all $n.$ By Beurling-Helson's Theorem, 
\[
|\mu(p)|=1, \mbox{ for some integer $p$, and} \mu(k)=0, \mbox{if $k\neq p$}. 
\]
Therefore, we have that 
\[
\mu=we^{(p)}, 
\]
for some integer $p$ and a complex number $w$ with modulus $1$ and that 
\[
h(e^{(1)})=we^{(p)},\quad h(e^{(j)})=[h(e^{(1)})]^{*j}=\mu^{*j}
=w^je^{(jp)}. 
\]
Finally, let $\rho\in l^1.$ Then, by the linearity and the continuity of $
\rho,$ 
\[
h(\rho)=h(\sum \rho(j)e^{(j)})=\sum \rho(j)\mu^{*j} =\sum
\rho(j)w^je^{(jp)}.
\]
Thus, every homomorphism of $l^1$ into $l^1$ has necessarily the above  form.

\vspace{0.1in}

{\bf Homomorphisms of $A(T)$ into $A(T)$}

\vspace{0.1in}

Let $h$ be a homomorphism of $A(T)$ into itself. Then one can show that $h$
necessarily has a trivial form, i.e., for any $g\in A(T),$ 
\[
h(g)(e^{it})=(g\circ m)(e^{it}), 
\]
where $m(e^{it})$ is the image of $e^{it}$ under $h.$

In fact, 
\[
h(g)=\sum \hat{g}(j)h(e^{ijt}) =\sum \hat{g}(j)(m(e^{it}))^j=g(m(e^{it}))
\quad (3). 
\]
That the series can be written as a composition of $g$ and $m$ can be easily
seen by replacing $e^{it}$ with $m(e^{it})$ in the Fourier series of $g$: $
g(e^{it})= \sum \hat{g}(j)e^{ijt}=\sum \hat{g}(j)(e^{it})^j.$

\begin{thm}
If $m(e^{it})\in A(T)$ is any mapping of $T$ into itself such that $
g\circ m\in A(T)$ for any $g\in A(T),$ then 
\[
m(e^{it})=we^{ipt} 
\]
for some integer $p$ and constant $w$ of modulus $1.$
\end{thm}

{\bf Proof:}\quad This is a restatement of Beurling- Helson's theorem. If $
g\circ m\in A(T),$ then we see from (3) that $|m^j(e^{it})|\leq K$ for $
j=0,\pm 1,\cdots.$ Let $\mathcal{F}(m)=\mu\in l^1.$ Then $\|\mu^{*j}\|_{l^1}
\leq K.$ Using Beurling-Helson's theorem, $\mu$ has a special form and so
does $m.$ 

\vspace{0.1in}

Let $m$ be any mapping of $T$ into itself. Define the mapping $h$ as 
\[
h(f)=f\circ m=\sum \hat{f}(j)(m(e^{it}))^j. 
\]
Obviously, $h$ is a homomorphism of $A(T)$ into $A(T).$ Define a
homomorphism of $l^1$ into $l^1$, denoted by $h^{\prime},$ in such a way
that 
\[
h^{\prime}(\mathcal{F}f)=\mathcal{F}(h(f)) 
\]
for all $f\in A(T).$ This can be written as 
\[
h^{\prime}(\rho)=\sum \rho(j)\mu^{*j} 
\]
for $\rho\in l^1,$ where $\mu=\mathcal{F}(m)\in l^1.$

\vspace{0.1in}

\begin{thm}
Assume that $m$ is a mapping of $T$ into $T$ such that $f\circ m\in A(T)$
whenever $f\in A(T).$ Then $h(f)=f\circ m$ is a homomorphism of $A(T)$ into $
A(T).$ Define $h^{\prime}$ as 
\[
h^{\prime}(\mathcal{F}f)=\mathcal{F}(h(f)),\qquad (4)
\]
that is, for $\rho\in l^1,$ 
\[
h^{\prime}(\rho)=\sum_{j}\rho(j)\mu^{*j}, 
\]
where $\mu=\mathcal{F}(m).$ Then $h^{\prime}$ is a continuous homomorphism of $
l^1$ into $l^1.$
\end{thm}

{\bf Proof:}\quad We show that $h^{\prime}$ is a homomorphism of $l^1$ into $
l^1.$ Let $\rho_k\in l^1$ and let $f_k=\mathcal{F}^{-1}(\rho_k),$ $k=1,2.$ Then 
\begin{eqnarray*}
h^{\prime}(\rho_1*\rho_2)=h^{\prime}(\mathcal{F}f_1*\mathcal{F}f_2)&=& h^{\prime}(
\mathcal{F}(f_1\cdot f_2) \\
=\mathcal{F}(h(f_1\cdot f_2)) &=&\mathcal{F}(h(f_1))*\mathcal{F}(h(f_2))=h^{\prime}(
\rho_1)*h^{\prime}(\rho_2).
\end{eqnarray*}

By the closed graph theorem, it is enough to show that if $
\rho^{(n)}\rightarrow \rho$ ($l^1(Z)$) and $h^{\prime}(\rho^{(n)})
\rightarrow \rho^{\prime}$ ($l^1(Z)$), then $\rho^{\prime}=h^{\prime}(\rho).$

\begin{eqnarray*}
\mathcal{F}^{-1}(\rho^{\prime})&=&\lim_{n\rightarrow\infty}\mathcal{F}^{-1}
(h^{\prime}(\rho^{(n)})) \\
&=&\lim_{n\rightarrow\infty}\mathcal{F}^{-1} (\sum \rho^{(n)}(j)\mu^{*j}) \\
&=&\lim_{n\rightarrow\infty}(\sum \rho^{(n)}(j)m(e^{it})^j) \\
&=&\sum \rho(j)m(e^{it})^j=\sum \rho(j)\mathcal{F}^{-1}(\mu^{*j}) \\
&=&\mathcal{F}^{-1}(\sum \rho(j)\mu^{*j}) =\mathcal{F}^{-1}(h^{\prime}(\rho)).
\end{eqnarray*}

The first equality holds because $h^{\prime}(\rho^{(n)})\rightarrow
\rho^{\prime}$ in $l^1(Z).$ The fourth equality holds because $
\rho^{(n)}\rightarrow \rho$ ($l^1(Z)$) and $|m(e^{it})|\leq 1.$

\begin{thm}
Let $h^{\prime}$ be any homomorphism of $l^1$ into $l^1.$ Then there exists $
m:T\rightarrow T$ so that for any $\rho\in l^1,$ 
\[
\mathcal{F}^{-1}[h^{\prime}(\rho)]=\mathcal{F}^{-1}(\rho)\circ m. 
\]
It follows that any homomorphism of $l^1$ into $l^1$ is of the form (4),
where $h(f)=f\circ m.$
\end{thm}

{\bf Proof:}\quad For each $q\in T,$ define 
\[
H_q(f)=[\mathcal{F}^{-1}(h^{\prime}(\rho))](q), 
\]
where $\rho=\mathcal{F}(f).$ Then $H_q$ is a homomorphism of $A(T)$ into $C.$
In fact, if we denote $\rho=\mathcal{F}f$ and $\eta=\mathcal{F}g,$ then 
\begin{eqnarray*}
H_q(f\cdot g)&=&[\mathcal{F}^{-1}(h^{\prime}(\rho*\eta))](q) \\
&=&[\mathcal{F}^{-1}(h^{\prime}(\rho)*h^{\prime}(\eta))](q) \\
&=&[\mathcal{F}^{-1}(h^{\prime}(\rho))\cdot \mathcal{F}^{-1}(h^{\prime}(\eta))](q)
\\
&=&H_q(f)\cdot H_q(g).
\end{eqnarray*}
Moreover, by the Gelfand theorem, $H_q$ must be continuous. Therefore, $H_q$
is an evaluation of the form $H_q(f)=f(q^{\prime})$ for some $q^{\prime}\in
T.$ Define $m: T\rightarrow T$ such that $q^{\prime}=m(q).$ Then, for any $
\rho\in l^1,$ let $f=\mathcal{F}^{-1}(\rho),$ we have 
\[
[\mathcal{F}^{-1}(h^{\prime}(\rho))](q)=H_q(f)=f(q^{\prime})=f\circ m(q) 
\]
for all $q\in T.$ 

\section*{Acknowledgements}
These notes were written by the first author in preparation for a series of talks given on harmonic analysis through a succession of seminars at the mathematics department of California State University Sacramento (CSUS). Later on, the second author joined in a full collaborative effort to revise and edit the notes and make them appropriate for publication as a graduate level textbook. We are distinctly grateful to the faculty of the mathematics department of CSUS for their helpful insights and support in preparation of these notes. In writing the present manuscript, we would also like to acknowledge that we were greatly inspired by Professor Henry Helson's classic book, Harmonic Analysis. Finally, we would like to express our sincere appreciation to Professor Calixto Calderon of the University of Illinois at Chicago for reviewing the final draft of the manuscript and making helpful suggestions to improve it.

\section*{Dedication}
The first author would like to sincerely express his gratitude to Mrs. Zhenyan Zhou, his late wife, for her affectionate support and encouragement during the writing of these notes. Without her unwavering and long time support, the present work would not have been possible. The second author also would like to acknowledge and express his gratefulness to Mrs. Mahin Aliabadi Siadat, his late mother, for her never-ending encouragement and loving support to persist in this collaboration, towards its successful conclusion. We dedicate the present work to these highly honorable and dedicated women. Although our loved ones are no longer with us, their memories will for ever last in our hearts.


\begin{thebibliography}{99}

\bibitem{1} R. L. Wheeden and A. Zygmund, Measure and Integral: An Introduction to Real Analysis,    	Marcel Dekker, Inc., New York and Basel, 1977.  
\bibitem{2} Helson, Harmonic Analysis (second edition), The Wadsworth and Brooks/Cole 	Mathematics Series, 1991.   	  
\bibitem{3} E. Hewitt and K. Stromberg, Real and Abstract Analysis, Springer-Verlag, Berlin, 	Heidelberg, New York, 1965.  
\bibitem{4} A. Zygmund, Trigonometric Series (second edition), Volume II, Cambridge University Press, 	Cambridge, New York, New Rochelle, Melbourne, Sydney, 1988. 
\bibitem{5} A. Zygmund, Trigonometric Series (second edition), Volume I, Cambridge University Press, 	Cambridge, New York, New Rochelle, Melbourne, Sydney, 1988. 


\end{thebibliography}
\end{document}